\documentclass[11pt,english,a4paper]{amsart} 
\pdfoutput=1

\usepackage[T1]{fontenc}
\usepackage[utf8]{inputenc}
\usepackage[a4paper]{geometry}
\usepackage{mathrsfs}
\usepackage{amsmath}
\usepackage{amssymb}
\usepackage{graphicx}
\usepackage{caption}
\usepackage{subcaption}
\usepackage{color}
\usepackage{a4wide}
\usepackage{amscd}
\usepackage{latexsym}
\usepackage{dsfont}
\usepackage{etoolbox} 
\usepackage{xparse} 
\usepackage{mathtools}
\usepackage{todonotes}
\usepackage{hyperref}
\usepackage{colonequals}
\usepackage{siunitx}

\makeatletter

\usepackage{fancyhdr}
\usepackage{esint}

\makeatother

\usepackage{babel}

\providecommand{\acknowledgementname}{Acknowledgement}
\providecommand{\definitionname}{Definition}
\providecommand{\examplename}{Example}
\providecommand{\lemmaname}{Lemma}
\providecommand{\remarkname}{Remark}
\providecommand{\theoremname}{Theorem}

\definecolor{grey}{rgb}{0.5,0.5,0.5}

\definecolor{mhcol}{rgb}{0,0.7,0}

\definecolor{jpcol}{rgb}{0.7,0.7,0}

\definecolor{lightblue}{rgb}{0.5,0.5,1}

\DeclareMathOperator{\Divergence}{div}

\newcommand{\op}[1]{\operatorname{#1}}
\newcommand{\st}{\;\vert\;}
\newcommand{\midst}{\;\middle\vert\;}

\definecolor{plot-green}{rgb}{0.4660, 0.6740, 0.1880}
\definecolor{plot-blue}{rgb}{0, 0.4470, 0.7410}
\definecolor{plot-orange}{rgb}{0.8500, 0.3250, 0.0980}
\definecolor{plot-red}{rgb}{0.6350, 0.0780, 0.1840}

\global\long\def\R{\mathbb{R}}
\global\long\def\N{\mathbb{N}}

\global\long\def\cS{\mathcal{S}}
\global\long\def\cT{\mathcal{T}}
\global\long\def\cN{\mathcal{N}}

\global\long\def\cC{\mathcal{C}}

\global\long\def\cE{\mathcal{E}}
\global\long\def\cJ{\mathcal{J}}

\global\long\def\cB{\mathcal{B}}

\global\long\def\cT{\mathcal{T}}
\global\long\def\cW{\mathcal{W}}
\global\long\def\cV{\mathcal{V}}

\newcommand{\abs}[1]{\left\vert #1 \right\vert}

\global\long\def\bodyCount{I}
\global\long\def\dim{d}
\global\long\def\Rd{\R^{\dim}}

\newcommand\closure[1]{\overline{#1}}
\newcommand{\frictionBoundary}[1][]{\Gamma^{F}_{#1}}
\newcommand{\dirichletBoundary}[1][]{\Gamma^{D}_{#1}}
\newcommand{\neumannBoundary}[1][]{\Gamma^{N}_{#1}}
\newcommand{\pen}[1]{\left[ #1 \right]}
\newcommand{\Ten}[1]{\boldsymbol{#1}}
\newcommand{\stressn}[1][]{\stressSymbol_n^{#1}}
\newcommand{\barstressn}[1][]{\bar{\stressSymbol}_n^{#1}}
\global\long\def\stressSymbol{\sigma}
\newcommand{\stresst}{\sigma_t}
\newcommand{\normal}{n}
\newcommand{\Vm}{V_m(\theta)}

\newcommand{\dieterich}{\psi_{\rm Dieterich}}
\newcommand{\ruina}{\psi_{\rm Ruina}}

\newcommand{\strain}{\Ten{\varepsilon}}
\newcommand{\stress}{\Ten{\sigma}}

\newcommand{\newtonCorr}{\delta u}
\newcommand{\Setminus}{\setminus}
\newcommand{\Id}{\operatorname{Id}}

\global\long\def\u{u}
\global\long\def\v{v}
\global\long\def\w{w}
\global\long\def\e{e}

\newcommand{\normalu}{\normal^{\u}}
\newcommand{\penu}[2][\u]{\pen{#2}^{#1}}
\newcommand{\frictionBoundaryU}[1][\u]{\Gamma^{F,#1}}
\newcommand{\cCu}{\cC^{\u}}
\newcommand{\Hu}[1][\u]{H_0^{#1}}
\newcommand{\Huold}{\Hu[\u_{n-1}]}
\newcommand{\Suold}{\cS_0^{\u_{n-1}}}

\newcommand{\piDisc}{\tilde{\pi}}
\newcommand{\penDisc}[1]{\widetilde{\pen{#1}}}
\newcommand{\frictionBoundaryDisc}{\tilde{\Gamma}^F}
\newcommand{\NodesNM}{\widetilde{\cN}^{F}}
\newcommand{\CDisc}{\tilde{\cC}}
\newcommand{\innerCDisc}[1]{\langle #1 \rangle_{\CDisc}}

\newcommand{\TOP}[1]{#1_{T}}
\newcommand{\TOPONE}[1]{#1_{T,1}}
\newcommand{\TOPI}[1]{#1_{T,\bodyCount-1}}
\newcommand{\TOPi}[1]{#1_{T,i}}

\newcommand{\BOT}[1]{#1_{B}}
\newcommand{\BOTONE}[1]{#1_{B,1}}
\newcommand{\BOTI}[1]{#1_{B,\bodyCount-1}}
\newcommand{\BOTi}[1]{#1_{B,i}}

\usepackage{amsthm}
\theoremstyle{remark}
\newtheorem*{acknowledgement*}{\protect\acknowledgementname}

\theoremstyle{plain}
\newtheorem{theorem}{\protect\theoremname}[section]

\theoremstyle{plain}

\theoremstyle{plain}

\theoremstyle{plain}

\theoremstyle{plain}
\newtheorem{problem}[theorem]{\protect Problem}

\theoremstyle{plain}
\newtheorem{lemma}[theorem]{\protect\lemmaname}

\theoremstyle{plain}
\newtheorem{proposition}[theorem]{\protect Proposition}

\theoremstyle{plain}

\theoremstyle{plain}

\begin{document}

\title{Numerical simulation of multiscale fault systems with rate- and state-dependent friction}

\author{Carsten Gr\"aser, Ralf Kornhuber and Joscha Podlesny}

\address{Carsten Gr\"aser, Institut f\"ur Mathematik, 
 Freie Universit\"at Berlin, 
 14195 Berlin, Germany,
 \href{https://orcid.org/0000-0003-4855-8655}{orcid.org/0000-0003-4855-8655}
}
\email{graeser@mi.fu-berlin.de}

\address{Ralf Kornhuber, Institut f\"ur Mathematik, 
 Freie Universit\"at Berlin, 
 14195 Berlin, Germany}
\email{kornhuber@math.fu-berlin.de}

\address{Joscha Podlesny, Institut f\"ur Mathematik, 
 Freie Universit\"at Berlin, 
 14195 Berlin, Germany}
\email{podlesjo@math.fu-berlin.de}


\date{}

\begin{abstract}
We consider the deformation of a geological structure with non-intersecting faults
that can be represented by a layered system of viscoelastic bodies 
satisfying rate- and state-depending friction conditions along the common interfaces.
We derive a mathematical model that
contains classical Dieterich- and Ruina-type friction as special cases
and accounts for possibly large tangential displacements.
Semi-discretization in time by a Newmark scheme leads to a coupled system of 
non-smooth, convex minimization problems for  rate and state to be solved in each time step.
Additional spatial discretization by a mortar method and piecewise constant finite elements
allows for the decoupling of rate and state by a fixed point iteration and 
efficient algebraic solution of the rate problem by truncated non-smooth Newton methods.
Numerical experiments with a spring slider and a layered multiscale system illustrate 
the behavior of our model as well as the efficiency and reliability of the numerical solver.
\end{abstract}

\maketitle

\thanks{
This research has been funded by Deutsche Forschungsgemeinschaft (DFG)
through grant CRC 1114 ''Scaling Cascades in Complex Systems'', Project Number 235221301,
Project B01  ''Fault networks and scaling properties of deformation accumulation''.
}

\section{Introduction}

%
Stress accumulation and release in geological fault networks play a crucial role in earthquake dynamics. 
The phenomenology of faults is ranging from subduction zones like the  Nasca plate and strike-slip faults like the San Andreas fault 
to multiscale fault systems like the Atacama zone. 
Strongly varying time scales between the occurrence and duration of  slip events
suggest to complement experimental studies in the field (or in the lab~\cite{rosenau2009shocks}) by numerical simulations. 

In the underlying mathematical description, the Dieterich–Ruina model of rate- and state-dependent friction (RSF)~\cite{Ruina83}
has become a standard for the frictional behaviour along the faults~\cite{ben2008collective,ranjith1999stability,rice2001rate}.
It can be regarded as an extension of simple Tresca friction with  rate- and state dependent friction coefficient $\mu = \mu(V, \theta)$ 
that is  increasing/decreasing with  increasing/decreasing sliding velocity or slip rate $V$ 
involving some relaxation effect as expressed by the state $\theta$. 
The variational structure of RSF has been identified and first exploited by Pipping et al.~\cite{PippingSanderKornhuber15} 

%
The simulation of rupture and slip events in seismic hazard analysis has quite a history (cf., e.g.,~\cite{barbot2012under,cochard1994dynamic,kato2004interaction,lapusta2000elastodynamic} and the references cited therein).
Utilizing a discontinuous Galerkin (DG) scheme in space in connection with  arbitrary high-order (ADER) time integration,
de la Puente et al~\cite{de2009dynamic} developed a numerical method for the dynamic simulation of  slip events.
This method was later generalized to three space dimensions~\cite{pelties2012three} and cast into the software package SeisSol
that was successfully utilized for the simulation, e.g., of the 2016 Kaik{\=o}ura earthquake cascade~\cite{ulrich2019dynamic}.
More recently, a different approach based on a diffuse representation of faults was first applied to subduction zones~\cite{herrendorfer2015earthquake,sobolev2017modeling,van2013seismic} and later extended to strike-slip faults~\cite{dal2018seismic}.
This approach has the potential to allow for much more complicated fault systems
because the faults have to be no longer resolved exactly by the underlying finite element mesh.
However, this advantage  currently comes with high computational cost
due to a lack of efficient algebraic solution techniques.

%
In this paper, we extend a variational approach to the simulation of subduction zones~\cite{pipping2016efficient}
to a  layered fault system with RSF.
More precisely, we consider the deformation of a geological structure with non-intersecting faults
that can be represented by $\bodyCount$ viscoelastic bodies
undergoing small viscoelastic deformations and large tangential displacements  with RSF contact conditions.
Assuming existence of a sufficiently regular contact mapping, we formulate a general mathematical model
that contains  Dieterich--Ruina friction as a special case.
Fault opening is forbidden for notational convenience, but could be included in a straightforward way.

Time discretization is performed by a classical Newmark scheme, 
resulting in a coupled system of non-smooth, convex minimization problems that has to be solved in each time step.
Decoupling this system by a fixed point iteration leads to a problem for the velocity (and thus for the rate)
with given state, and an independent state problem with given rate.
Both the rate and the state problem can be rewritten as convex minimization problems admitting unique solutions.
For a related  coupled problem with unilateral contact,
as arising from the mathematical description of subduction zones,
existence  and uniqueness of solutions was  established by Pipping et al.~\cite{Pipping2014,PippingSanderKornhuber15} 
using  fixed point arguments.

Spatial discretization of the rate problem is performed by a mortar method in the spirit of Krause and Wohlmuth~\cite{krause2002dirichlet,wohlmuth2003monotone,wohlmuth2011variationally}.
This approach has the advantage that it provides nodal block separation of  the non-smooth nonlinearity
which allows for direct application of globally convergent {Truncated Non-smooth Newton Multi-Grid} (TNNMG) methods~\cite{Graeser2011,GraeserSackSander2009,GraeserSander2019}.
The  state problem is discretized by  piece-wise constant finite elements.
For given rate, the  resulting algebraic problem decouples into independent scalar problems for each of the nodal values,
which can be solved, e.g., by simple bisection or even explicitly. 

In our numerical experiments, we  consider a  spring slider with  $\bodyCount = 2$ bodies 
and  a layered network  with $\bodyCount = 5$ bodies separated by 4 faults subject to prescribed velocities at the upper boundary.
We perform self-adaptive time stepping to efficiently resolve strongly varying velocities during  loading, rupture, and sliding.
Spatial discretization is based on triangulations as obtained by adaptive refinement concentrated at the faults.
The associated hierarchy of finite element spaces is used for the algebraic TNNMG solver of the rate problems with given state
as arising in the fixed point iteration mentioned above.

For the spring slider we observe the periodic occurrence of mostly unilateral slip events,
similar to related simulations of subduction zones~\cite{pipping2016efficient}. 
These slip events are nicely captured by adaptive time stepping, while the number of outer fixed point iterations
and inner multigrid iteration remains almost the same for all time steps. 
Simulation of the layered network  exhibits an interesting coincidence of periodic slip events along the upper fault 
with loading phases and oscillatory behavior on the others.
We observe essentially the same efficiency of time stepping, fixed point iteration, and multigrid  as for the spring slider
which illustrates the robustness of our numerical solution procedure, also with respect to the number of faults.


\section{Mathematical modelling}
\subsection{A layered fault system with rate-and state-dependent friction}
We consider a  geological structure containing a system of  faults
which is represented by a deformable body with reference domain $\Omega \subset \Rd$,  $\dim = 2,3$,
that, along the faults, is decomposed  into $\bodyCount$ subdomains $\Omega_i$, $i=1, \dots,\bodyCount$,
\[
\overline{ \Omega} = \bigcup_{i=1}^{\bodyCount} \overline\Omega_i .
\]
We  assume that these subdomains are
 non-empty, bounded Lipschitz domains,  do not penetrate each other and are layered 
in the sense that at most two subdomains are in contact at any point in $\Rd$ (see Figure~\ref{fig:layeredBlocks}). 
Then, the subdomains can be ordered such that there is a 
common interface 
$\frictionBoundary[{i,i+1}] = \closure{\Omega}_i \cap \closure{\Omega}_{i+1}$,  $i=1,\dots, \bodyCount -1$, and 
all other intersections of subdomains are empty.
Setting $\frictionBoundary[{0,1}]=  \frictionBoundary[{\bodyCount,\bodyCount +1}] =\emptyset$ for notational convenience,
the boundary $\partial \Omega_i$ of  $\Omega_i$ is disjointly decomposed according to   
$\partial \Omega_i = \dirichletBoundary[i] \cup \neumannBoundary[i] \cup (\frictionBoundary[{i-1,i}]\cup  \frictionBoundary[{i,i+1}])$, 
into a Dirichlet, a Neumann, and a contact  boundary, respectively.
We set
\begin{equation*}
    \dirichletBoundary = \bigcup_{i=1}^{\bodyCount} \dirichletBoundary[i], \quad
    \neumannBoundary = \bigcup_{i=1}^{\bodyCount} \neumannBoundary[i], \quad
    \frictionBoundary = \bigcup_{i=1}^{\bodyCount-1} \frictionBoundary[{i,i+1}].
\end{equation*}

For $\v=(v_1, \dots, v_\bodyCount)$ with $v_i: \Omega_i \to \Rd$, $i=1,\dots,\bodyCount$, we define the restrictions
$\TOP{\v} =  (\TOPONE{v}, \dots, $ $\TOPI{v})$ and $\BOT{\v}=  (\BOTONE{v}, \dots, \BOTI{v})$ 
of  $\v$ to $\frictionBoundary$ with  
\begin{align*}
  \TOPi{v}= v_{i+1}|_{\frictionBoundary[i,i+1]}, \quad
  \BOTi{v}= v_i|_{\frictionBoundary[i,i+1]} \qquad
 i =1, \dots, \bodyCount-1,
\end{align*}
denoting the restrictions from the top $\Omega_{i+1}$ and the bottom $\Omega_i$, respectively (see Figure~\ref{fig:layeredBlocks}).
It is convenient to identify  $\TOP{\v} =  (\TOPONE{v}, \dots, $ $\TOPI{v})$ and $\BOT{\v}=  (\BOTONE{v}, \dots, \BOTI{v})$ 
with  functions $\TOP{\v}$ and $\BOT{\v}$ defined on $\frictionBoundary$
by  $\TOP{\v}|_{\frictionBoundary[i,i+1]}  = \TOPi{v}$ and $\BOT{\v}|_{\frictionBoundary[i,i+1]}  = \BOTi{v}$,  $i=1,\dots,\bodyCount-1$.
Let $\normal =(\normal_1, \dots, \normal_\bodyCount)$ where $n_i \in \Rd$  stands for the
outer normal to $\Omega_i$, $i=1,\dots \bodyCount$. 
Note that $\normal_i$ is an inward normal to $\Omega_{i+1}$ on $\frictionBoundary[{i,i+1}]$.
In particular,  $\TOP{\normal} = (\TOPONE{\normal}, \dots, \TOPI{\normal})$ and  $\BOT{\normal} = - \TOP{\normal}$
are  top and bottom normals on $\frictionBoundary$, respectively.
In the following,  most quantities will be  defined in terms of the bottom side.

\begin{figure}
	\centering 
    \input{gfx/layeredblocks.pgf}
    \caption{A  fault system $\frictionBoundary$ with $\bodyCount=4$ subdomains and $\bodyCount - 1 = 3$ layered faults}
    \label{fig:layeredBlocks}
\end{figure}

We suppose that a body force $f$ acts on all of  $\Omega$ and  surface forces $f^N$ act on the  Neumann boundary $\neumannBoundary$. 
On the Dirichlet boundary $\dirichletBoundary$ the velocity $\dot u(t)$ of the displacement field $u(t)$  
of the deformable body $\Omega$ is fixed at all time instants $t>0$.
We set  $u(t)= \dot u(t)=0$ on  $\dirichletBoundary$ for convenience,
though all further considerations can be  generalized to the inhomogeneous case in a straightforward way.

We assume that the boundary forces are compressive in the sense that no fault opening occurs.
This means that  neighboring bodies $\Omega_i$ and $\Omega_{i+1} $, $i=1, \dots,\bodyCount -1$, 
remain in contact throughout the evolution.

  We consider a deformation field
  \begin{align*}
    \u = (u_1, \dots,u_{\bodyCount}) \in  H^1(\Omega_1)^d\times \cdots \times H^1(\Omega_{\bodyCount})^d 
  \end{align*}
  where $\u_i$ is the deformation of the subdomain $\Omega_i$.
  Throughout the paper we assume that the deformations $\u_i$
  within each subdomain $\Omega_i$ are small, while the relative
  displacement between different subdomains can be large.
  Thus we will use a (geometrically) linear elastic approach
  inside of the subdomains $\Omega_i$ while the coupling conditions
  have to take care of potentially large deformations.

  Large deformation coupling conditions will be defined in terms of the
  deformed subdomains.
  Given the deformation fields $\u_i$ the associated displacements
  are given by $\Id + \u_i$ leading to the deformed subdomains $(\Id + \u_i)(\Omega_i)$.
  The actual contact surface of the deformed subdomains is then given by
  $\cCu = (\Id + \BOT{\u})(\frictionBoundary) \cap (\Id + \TOP{\u})(\frictionBoundary)$.
  In the following, we assume that each $\Id + \u_i$ is injective, i.e. that each
  $\u_i$ is regular enough to avoid self-penetration of $\Omega_i$.
  Furthermore, we assume that deformations are small, such that
  different surfaces $\frictionBoundary[i,i+1]$ do not
  get in contact after deformation.
  Then, the deformed contact boundary can be pulled back to the bottom
  and top reference domain according to
  \begin{align*}
    \BOT{\frictionBoundaryU} &= (\Id + \BOT{\u})^{-1}(\cCu) \subset \frictionBoundary,&
    \TOP{\frictionBoundaryU} &= (\Id + \TOP{\u})^{-1}(\cCu) \subset \frictionBoundary.
  \end{align*}
  In the following, we will parameterize the top reference domain $\TOP{\frictionBoundaryU}$
  over the bottom one $\BOT{\frictionBoundaryU}$ by the bijective contact mapping
  \begin{align*}
    \pi^{\u}: \BOT{\frictionBoundaryU} \to \TOP{\frictionBoundaryU},
    \qquad
    \pi^{\u} = (\Id+\TOP{\u})^{-1} \circ (\Id+ \BOT{\u}),
  \end{align*}
  which maps each bottom point $x \in \BOT{\frictionBoundaryU}$
  to the unique top point $y \in \TOP{\frictionBoundaryU}$, such that the corresponding
  displaced points $(\Id + \BOT{\u})(x)$ and $(\Id + \TOP{\u})(y)$
  coincide.
  As a consequence, the deformed contact boundary $\cCu$
  can be parametrized  both over $\BOT{\frictionBoundaryU}$ using $\Id + \BOT{\u} = (\Id + \TOP{\u}) \circ \pi^\u$
  and over $\TOP{\frictionBoundaryU}$ using $\Id + \TOP{\u} = (\Id + \BOT{\u}) \circ (\pi^\u)^{-1}$.
  
  Now, consider any piecewise defined scalar or vector field
  \begin{align*}
    \v = (\v_1,\dots,\v_I)\in H^1(\Omega_1)^k\times \cdots \times H^1(\Omega_{\bodyCount})^k
  \end{align*}
  with $k=1$ or $k=d$.
  Then, we define the jump of $\v$ across the deformed contact boundary $\cCu$ 
  on the contact reference domain $\BOT{\frictionBoundaryU}$ according to
  \begin{align} \label{eq:CONTJUMP}
    \penu{\v}
      = \BOT{\v} - \TOP{\v} \circ \pi^{\u} \qquad \text{on }\BOT{\frictionBoundaryU}.
  \end{align}

  Contact conditions and friction laws will be stated in terms of
  normal and tangential components on the deformed contact boundary $\cCu$. 
  To this end, let $i = 1, \dots \bodyCount$ and denote by  $\normalu(x)$
  an outer normal to
  $(\Id + \u_i)(\Omega_i)$ at the point $(\Id + \u_i)(x)$, $x \in \frictionBoundary[i,i+1]\cap \BOT{\frictionBoundaryU}$,
  i.e. $\normalu$ is the pullback of an oriented normal field of
  the deformed contact boundary $\cCu = (\Id + \BOT{\u})(\BOT{\frictionBoundaryU})$
  to $\BOT{\frictionBoundaryU}$ using the map $\Id + \BOT{\u}$.
  Then we can decompose any
  vector field on $\BOT{\frictionBoundaryU}$
  according to its normal and tangential components with respect
  to the deformed configuration
  as
  \begin{align*}
    \v = \v_t + \v_\normal \normalu,
    \qquad \v_\normal = \v\cdot \normalu,
    \qquad \v_t = \v - \v_\normal \normalu.
  \end{align*}
  It is important to note that the tangential and normal component
  are defined in terms of the $\u$-dependent normal field $\normalu$
  that is defined (piecewise) on the deformed bottom subdomains $(\Id+\u_i)(\Omega_i)$
  and not with respect to the normal field of the reference subdomains $\Omega_i$.

  We state a closed-fault condition (no penetration  and  no  fault opening) by prescribing
  that the relative motion of the deformed subdomains $\Omega_i$
  and $\Omega_{i+1}$ is tangential to the actual contact surface $\cCu$, i.e.
  \begin{equation} \label{eq:NOFO}
    0 = \penu{\dot \u} \cdot \normal^{\u}.
  \end{equation}
  As a con\-se\-quence, the jump of  the relative tangential velocity satisfies
  \[
    \penu{\dot \u }_t = \penu{ \dot \u} - (\penu{\dot \u}\cdot\normalu) \normalu =  \penu{ \dot \u }\quad \text{on } \BOT{\frictionBoundaryU}.
  \]
  
  The closed-fault condition is complemented by the balance of normal forces
  \begin{align}
    \label{eq:normalForceBalanceRaw}
    \BOT{(\stress(\u)\normal)}  = -\omega^{\u} \TOP{(\stress(\u) \normal)} \circ \pi^{\u}
  \end{align}
  on $\BOT{\frictionBoundaryU}$,
  where $\stress(\u)$ denotes the stress tensor on $\Omega$.
  Note that the normal force
  $\TOP{(\stress(\u)\normal)}$ is a force per surface area,
  such that the change of the area element induced by the pullback to
  $\BOT{\frictionBoundaryU}$ using $\pi^{\u}$
  has to be compensated by the weighting factor
  $\omega^{\u} = \sqrt{\op{det}((D\pi^{\u})^TD\pi^{\u} + \BOT{\normal} \otimes \BOT{\normal})}$,
  while the minus sign compensates for the change of the normal
  direction on opposing sides.
  Note that the balance of normal forces \eqref{eq:normalForceBalanceRaw} can alternatively be phrased
  as a jump condition with a transformed weighting factor  $\omega^{\u}\circ (\pi^{\u})^{-1}$.

Utilizing  $\penu{\dot \u }_t =  \penu{\dot \u}$, we prescribe a rate- and state-dependent friction law of the form 
\begin{equation}\label{eq:frictionLaw}
  - \stresst \in  \partial_{\penu{\dot \u}} \phi(\penu{\dot{\u}}, \alpha) \quad \text{on } \BOT{\frictionBoundaryU}.
\end{equation}
Here, we used the decomposition $\BOT{(\stress(\u)\normal)} = \stresst + \stressn \normalu$
of the stress field $\BOT{(\stress(\u) \normal)}$
on the bottom side into its normal and tangential components
\begin{align*}
  \stressn &= \stressn(\u) = \BOT{(\stress(\u)\normal)} \cdot \normalu, &
  \stresst &= \stresst(\u) = \BOT{(\stress(\u)\normal)} - \stressn(\u)\normalu,
\end{align*}
respectively, and $\partial_{\penu{\dot \u}} \phi(\penu{\dot \u},\alpha)$ 
denotes the subdifferential of a state-dependent convex functional $\phi(\cdot, \alpha)$ to be described below.
Note that the stress vector $\BOT{(\stress(\u)\normal)}$ is computed with respect to the reference normal $\BOT{\normal}$,
while its decomposition in tangential and normal components is computed with respect to the deformed configuration
with the corresponding normal $\normalu$.
This reflects the fact that we assume small deformations within the subdomains
while the relative deformations of subdomains can be large.

For given relative slip rate $\abs{\penu{\dot \u}}$, the evolution of the state $\alpha$  is given by
\begin{equation}
- \dot \alpha = \partial_{\alpha} \psi(\alpha, \abs{\penu{\dot{\u}}}) \quad \text{on }\BOT{\frictionBoundaryU}, \qquad 
  - \dot \alpha = 0  \quad \text{on }  \frictionBoundary \Setminus \BOT{\frictionBoundaryU}
\end{equation}
with a second convex functional $\psi(\cdot,  \abs{\penu{\dot{\u}}})$.
Note that the state $\alpha$ remains constant on $ \frictionBoundary \Setminus\BOT{\frictionBoundaryU}$ where no contact occurs.

Assuming a visco--elastic Kelvin--Voigt material law, and fixing some final time $T_0>0$,
we are now ready to state the following formal description of the 
deformation of a body $\Omega$  with a layered fault system $\frictionBoundary$ and rate- and state-dependent friction.
\begin{problem}[Layered fault system with rate- and state dependent friction] \label{prob:CLASSICAL}
 Find 
  \[
  u \colon \; \Omega \times [0,T_0] \to \Rd \quad \text{and} \quad  \alpha \colon \: \frictionBoundary \times [0,T_0] \to \R
  \]
 such that
  \begin{align}
    \label{eq:kelvin-voigt} 
    \stress(u) & =  \Ten{A} \strain(\dot \u) + \Ten{B} \strain(\u)&& \text{in $\Omega\Setminus \frictionBoundary$} && \text{(Kelvin--Voigt material)} \\
    \label{eq:balance-of-momentum}
     \Divergence \stress(u) + f  & = \rho \ddot u && \text{in $\Omega\Setminus \frictionBoundary$} && \text{(balance of momentum)}\\
    \shortintertext{with boundary conditions,}
    \nonumber \u = \dot \u & = 0 && \text{on $\dirichletBoundary$}  && \text{(Dirichlet condition)}\\
    \nonumber \stress(\u) \normal & = f^N && \text{on $\Gamma^N$}  && \text{(Neumann condition)}\\
    \shortintertext{ frictional contact conditions,}
    \label{eq:NOPEN}\penu{\u} \cdot \normalu & = 0 && \text{on } \BOT{\frictionBoundaryU} && \text{(closed-fault condition)}\\
    \label{eq:normalForceBalance} 
     \BOT{(\stress(\u)\normal)}  &= -\omega^{\u} \TOP{(\stress(\u) \normal)} \circ \pi^{\u} && \text{on } \BOT{\frictionBoundaryU} && \text{(balance of normal forces)}\\
    \label{eq:1} -\stresst &\in \partial_{\penu{\dot \u}} \phi(\penu{\dot{\u}}, \alpha) &&\text{on } \BOT{\frictionBoundaryU} && \text{(state-dependent friction law)}\\
    \shortintertext{contact state condition,}
    \label{eq:2}
    - \dot \alpha &\in \partial_\alpha \psi(\alpha,\abs{\penu{\dot{\u}}}), &&\text{on } \BOT{\frictionBoundaryU}  && \text{(rate-dependent state law) }\\
    \shortintertext{and non-contact interface conditions}
       \label{eq:NONCONTACT_STATE} - \dot \alpha &= 0 && \text{on } \frictionBoundary \Setminus \BOT{\frictionBoundaryU} && \text{(non-contact state condition)}\\
    \label{eq:NONCONTACT_BOT} \BOT{(\stress(\u)\normal)} & = 0 && \text{on } \frictionBoundary \Setminus \BOT{\frictionBoundaryU} && \text{(bottom Neumann condition)}\\
    \label{eq:NONCONTACT_TOP} \TOP{(\stress(\u)\normal)} & = 0 && \text{on } \frictionBoundary \Setminus \TOP{\frictionBoundaryU} && \text{(top Neumann condition)}
  \end{align}
  holds for all $t \in [0,T_0]$. Here,  $\rho>0$ is a constant material density, 
  $\Ten{A}$ and $\Ten{B}$ stand for the viscosity and elasticity tensor, respectively,  
  and $\strain(v) = \frac{1}{2}(\nabla v + (\nabla v)^T)$ is the linearized strain or strain rate tensor. 
  In addition, we impose initial conditions on the displacement $u$, velocity $\dot u$, and state~$\alpha$.
\end{problem}

Throughout the following, we assume that the tensor fields $\Ten{A}$ und $\Ten{B}$
have the symmetry properties
\begin{align*}
  \Ten{A}_{ijkl} &= \Ten{A}_{klij}, &
  \Ten{A}_{ijkl} &= \Ten{A}_{jikl}, &
  \Ten{B}_{ijkl} &= \Ten{B}_{klij}, &
  \Ten{B}_{ijkl} &= \Ten{B}_{jikl}
\end{align*}
such that the stress tensor $\stress(\u)$ and the bilinear forms induced by
$\Ten{A}$ and $\Ten{B}$ are symmetric.
%

Note that Problem~\ref{prob:CLASSICAL} provides an extension of the model presented in~\cite{PippingSanderKornhuber15}
that describes unilateral frictional  contact  of  a deformable body  with a rigid foundation.
The tangential velocity relative to the fixed rigid foundation appearing in~\cite{PippingSanderKornhuber15} 
is  now replaced by the relative tangential velocity of adjacent deformable bodies.

A further extension to fault opening  can be performed by 
replacing \eqref{eq:NOPEN} by  the non-penetration condition $\penu{\u} \cdot \normalu \leq 0$
together with dynamical freezing and thawing of rate- and state-dependent friction \eqref{eq:1}, \eqref{eq:2} in case of opening or closing  faults.

For ease of notation we will skip the superscript and mostly write  $\pen{\cdot} = \penu{\cdot}$  in the sequel.


\subsection{The Dieterich--Ruina model}
The current form of the Dieterich--Ruina model of rate- and state-dependent friction goes back to~\cite{Ruina83} (see also \cite{doi:10.1029/2007JB005082,doi:10.1029/94GL01599,doi:10.1144/SP359.1,doi:10.1146/annurev.earth.26.1.643,doi:10.1073/pnas.93.9.3811,doi:10.1038/34097} and the papers cited therein).
It is based on the following ansatz for the friction coefficient 
\begin{equation} \label{eq:FRICCOEFF}
\mu^{*} (V, \theta) = \mu_0 + a \log\left( \frac{V}{V_0}\right) + b \log \left( \frac{V_0 \theta}{L}  \right)
\end{equation}
 that depends on the rate $V = |\pen{\dot u}|$ and involves positive parameters $\mu_0$, $V_0$ $a$, $b$, and $L$.
 
 It is complemented by a suitable evolution of the  state $\theta > 0$. 
 Here, most popular choices are
 \begin{align}
  \label{eq:DIETERICH} \dot \theta  &= 1 - \frac VL\theta &&\text{(Dieterich's law)}\\
 \shortintertext{and} 
 \label{eq:RUINA}  \dot \theta &= - \frac VL\theta \log\left(\frac{V}{L}\theta\right)  &&\text{(Ruina's law).}
\end{align}
Following \cite{PippingSanderKornhuber15}, we briefly sketch how this completely phenomenological friction model translates 
into  a corresponding state-dependent friction law \eqref{eq:1} and a rate-dependent state evolution \eqref{eq:2} as postulated above.
Starting from collinearity of relative tangential velocity and stress
\[
 -  \stresst  |\pen{\dot u}|= \pen{\dot u}  |\stresst|, 
\]
we postulate the state equation $ |\stresst| = \mu^*(V,\theta) |\stressn|$ 
with normal stress $\stressn = \stressn(\u)$ to obtain
\begin{equation} \label{eq:DRCLASSIC}
-\stresst  = |\stressn|\mu^*(V,\theta) \frac{\pen{\dot u}}{|\pen{\dot u}|}.
\end{equation}
In analogy to Tresca friction, we now replace the solution dependent normal stress $|\stressn| = |\stressn(\u)|$ 
by a given parameter $|\barstressn|$~\cite{PippingSanderKornhuber15}.

As  $\mu^*(V,\theta)$ becomes negative and thus meaningless for
\[
0\leq V < \Vm = V_0 \exp\left( - \frac{\mu_0 + b + \log(\theta V_0 /L)}{a} \right)
\]
we replace  $\mu^*(V,\theta)$ by its regularization 
\[
\mu(V,\theta) = 
\left\{
\begin{array}{cc}
 \mu^*(V,\theta) &\text{ if } V \geq \Vm\\
 0 &\text{ otherwise }
 \end{array}
 \right .  .
\]
Then elementary calculations show that \eqref{eq:DRCLASSIC} takes the form \eqref{eq:1}  with the convex functional
$\phi^{*} (\cdot,\theta): \Rd \to \R$ defined by
\begin{equation} \label{eq:RATEPHI}
\phi^*(v,\theta) = 
\left\{
\begin{array}{cc}
  a |\barstressn| \left( |v| \log(|v| / \Vm ) -|v| + \Vm\right) &\text{ if } |v| \geq \Vm\\
 0 &\text{ otherwise }
 \end{array}
 \right .  .
\end{equation}

It remains to  show that the rate evolutions \eqref{eq:DIETERICH}  and \eqref{eq:RUINA} can be rewritten according to \eqref{eq:2}. 
Introducing the transformed state $\alpha = \log \theta$,  Dieterich's law  \eqref{eq:DIETERICH}  
takes the form \eqref{eq:2} with the scalar  convex function 
\begin{equation}\label{eq:DIETERICHLAW}
  \dieterich(\alpha, V)  = \frac{V}{L}\alpha + e^{-\alpha} .
\end{equation}
Ruina's law \eqref{eq:RUINA}  is recovered  in terms of  \eqref{eq:2} by the same transformation and the scalar  convex function 
\begin{equation}\label{eq:RUINALAW}
 \ruina(\alpha, V)    = \frac {V}{L} \left( \frac {1}{2} \alpha^2 + \log(V/L) \alpha \right).
 \end{equation}
 Inserting the transformed state $\theta = e ^{\alpha}$ into \eqref{eq:RATEPHI}, we obtain the corresponding rate functional
 \begin{equation} \label{eq:DRF}
 \phi(\cdot, \alpha) = \phi^*(\cdot, e^{\alpha}). 
 \end{equation}

\subsection{Weak formulation}
We consider the Hilbert space $H =  H^1(\Omega_1)^d \times \cdots \times H^1(\Omega_{\bodyCount})^d$ with the canonical inner product
$(\v,\w)_H = \sum_{i=1}^{\bodyCount} (\v_i,\w_i)_{H^1(\Omega_i)^d}$, $\v_i, \w_i \in H^1(\Omega_i)$, $i=1,\dots, \bodyCount$,
and introduce the closed linear subspace
\begin{equation*}
H_0 = \{ \v \in H \;|\; 
\v = 0 \text{ on }\dirichletBoundary\}
\end{equation*}
of admissible displacements
respecting the Dirichlet boundary conditions.
The normal jump condition is incorporated in the closed
affine subspace
\begin{align} \label{eq:SOLSPACE}
  \Hu = \{v \in H_0 \st \penu{v} \cdot \normalu = 0\}.
\end{align} 
With the tensors $\Ten{A}$, $\Ten{B}$ taken from \eqref{eq:kelvin-voigt}, we introduce the bilinear forms
\begin{equation} \label{eq:AB}
a(\v, \w) = \int_{\Omega \Setminus \frictionBoundary}  \Ten{A} \strain( \v) \colon \strain(\w)\; dx, \quad
b(\v, \w) = \int_{\Omega \Setminus \frictionBoundary} \Ten{B} \strain( \v) \colon \strain(\w)\; dx, \qquad v, w \in H_0,
\end{equation}
involving the linear strain tensor $\strain(\v)= \frac{1}{2}(\nabla v + (\nabla v)^T)$, together with the linear functional
\begin{equation} \label{eq:ELL}
\ell(v) = \int_{\Omega} f v \; dx + \int_{\neumannBoundary} f^N  v \; ds, \qquad v \in H_0.
\end{equation}
To ensure that the bilinear forms are well defined,
we assume that the tensor fields $\Ten{A}$ und $\Ten{B}$
are uniformly elliptic in the sense that the bilinear forms
induced by $\Ten{A}(x)$ and $\Ten{B}(x)$ on the space of symmetric
$d \times d$ matrices are elliptic with constants independent of $x \in \Omega$.

By inserting the stress--strain relation \eqref{eq:kelvin-voigt} into the balance of momentum \eqref{eq:balance-of-momentum},
testing   with $v - \dot u$, integrating by parts, and exploiting the symmetry of $\stress(u)$
together with the boundary conditions  on $\dirichletBoundary$ and $\neumannBoundary$
we then formally obtain 
\begin{multline} \label{eq:WEMO}
 \langle \rho \ddot u, v - \dot u \rangle  + a( \dot u, v - \dot u)  +  b(u, v - \dot u)
  - \ell(v -\dot u) \\
  =
    \int_{\frictionBoundary} \BOT{(\stress(u) \normal)} \cdot \BOT{(v - \dot u)} \; ds +
    \int_{\frictionBoundary} \TOP{(\stress(u) \normal)} \cdot \TOP{(v - \dot u)} \; ds
\end{multline}
for all $ v \in H_0$ and $t \in (0,T_0)$.  Here $\langle \cdot, \cdot \rangle$ stands for the pairing of $H_0$ with its dual  $H_0^*$.
Using the boundary conditions \eqref{eq:NONCONTACT_BOT} and \eqref{eq:NONCONTACT_TOP},
integral transformation from $\TOP{\frictionBoundaryU}$ to $\BOT{\frictionBoundaryU}$,
and the normal force balance \eqref{eq:normalForceBalance} we can rewrite the right hand side in \eqref{eq:WEMO} as
\begin{align*}
 & \int_{\BOT{\frictionBoundaryU}}
    \BOT{(\stress(u) \normal)} \cdot \BOT{(v - \dot u)} \; ds +
  \int_{\TOP{\frictionBoundaryU}}
    \TOP{(\stress(u) \normal)} \cdot \TOP{(v - \dot u)} \; ds\\
  =&
  \int_{\BOT{\frictionBoundaryU}}
    \BOT{(\stress(u) \normal)} \cdot
      \Bigl(\penu{v - \dot u}  + \TOP{(v - \dot u)} \circ \pi^\u\Bigr) \; ds +
  \int_{\BOT{\frictionBoundaryU}}
    \omega^{\u} \Bigl(\TOP{(\stress(u) \normal)} \cdot \TOP{(v - \dot u)}\Bigr)\circ \pi^\u \; ds\\
  =&
  \int_{\BOT{\frictionBoundaryU}}
    \BOT{(\stress(u) \normal)} \cdot \penu{v - \dot u} \; ds.
\end{align*}

For given state $\alpha$, we introduce the convex functional $\Phi^\u$ on $H_0$ according to
\begin{equation} \label{eq:RATEFUN}
  \Phi^\u ( \cdot ,  \alpha) = \int_{\BOT{\frictionBoundaryU}} \phi (\penu{\cdot},  \alpha) \; ds 
\end{equation}
with the convex functional $\phi$ taken from the friction law \eqref{eq:1}.
Now let $\v$ satisfy the closed-fault condition \eqref{eq:NOPEN}, i.e. $\v \in \Hu\subset H_0$.
Then, utilizing the  decomposition
$\BOT{(\stress(\u)\normal)} = \stresst + \stressn \normalu$
together with the friction law \eqref{eq:1}, the closed-fault condition
$\penu{v - \dot u} \cdot  \normalu=0$ on $\BOT{\frictionBoundaryU}$,
and the definition of subdifferentials,
we find that
\begin{equation} \label{eq:BINCL}
  \int_{\BOT{\frictionBoundaryU}} \BOT{(\stress(\u) \normal)} \cdot \penu{\v - \dot \u} \; ds
  \geq   \Phi^\u( \dot \u,  \alpha) -  \Phi^\u( \v , \alpha) .
\end{equation}

Now we insert   \eqref{eq:BINCL} into  \eqref{eq:WEMO}, in order to obtain the desired weak form of the rate equation
 \begin{equation} \label{eq:WEAKRATE}
 \langle \rho \ddot \u, \v - \dot \u \rangle  + a( \dot \u, \v - \dot \u)  + b(\u, \v - \dot \u)
   + \Phi^\u(\v, \alpha) - \Phi^\u(\dot \u, \alpha)  \geq  \ell(\v -\dot \u) \quad  \forall  v \in \Hu.
 \end{equation}
 
Similarly, for given velocity  $\dot u \in \Hu$ and thus given rate $\abs{\penu{\dot \u}}$,
we define the convex functional  $\Psi^\u$ on $L^2(\frictionBoundary)$ by 
\begin{equation} \label{eq:STATEFUN}
  \Psi^\u ( \cdot,  \dot u) = \int_{\BOT{\frictionBoundaryU}} \psi (\cdot, \abs{\penu{\dot \u}} ) \; ds 
\end{equation}
with the convex functional $\psi$ taken from the state law \eqref{eq:2},
and test  the state evolution \eqref{eq:2} with $\beta \in L^2 (\frictionBoundary)$ to obtain the weak  formulation
\begin{equation} \label{eq:WST}
  ( \dot \alpha, \beta - \alpha)_{L^2 (\frictionBoundary)} + \Psi^\u (\beta,  \dot u) - \Psi^\u (\alpha, \dot u) \geq  0
		\qquad \forall \beta \in L^2 (\frictionBoundary).
\end{equation}
This formulation automatically satisfies the non-evolution condition \eqref{eq:NONCONTACT_STATE}
for the state $\alpha$ on the non-contact boundary $\frictionBoundary \Setminus \frictionBoundaryU$
since $\Psi^\u$ is defined on $L^2(\frictionBoundary)$ but only
depends on values of $\alpha$ on the contact boundary $\frictionBoundaryU$.

We are now ready to state the weak formulation of Problem~\ref{prob:CLASSICAL}
\begin{problem}[Weak formulation] \label{eq:WEAK}
Find 
\[
\u \in H^1((0,T_0), H_0) \cap H^2((0,T_0), H_0^*) \text{ and }\alpha \in H^1((0,T_0), L^2(\frictionBoundary))
\]
such that
$\dot \u \in \Hu$ and
\begin{align}
    \label{eq:WEAKRATE_COUPLED} 
\langle \rho \ddot \u, \v - \dot \u \rangle  + a(\dot \u, \v - \dot \u)  + b(\u, \v - \dot \u) 
  + \Phi^\u (\v, \alpha) - \Phi^\u (\dot \u, \alpha) & \geq  \ell(\v -\dot \u) \quad & \forall &v\in \Hu,\\
 \label{eq:WEAKSTATE}
  ( \dot \alpha, \beta - \alpha)_{L^2 (\frictionBoundary)} + \Psi^\u (\beta, \dot \u ) - \Psi^\u (\alpha,   \dot \u ) & \geq  0
		\quad & \forall & \beta \in L^2 (\frictionBoundary)
\end{align}
holds for almost all $t \in (0,T_0)$ together with initial conditions
\begin{equation} \label{eq:INIT}
u(0)=u_0, \quad \dot u(0)= \dot u_0,  \quad \alpha(0)= \alpha_0
\end{equation}
with given $u_0 \in \Hu[\u_0]$ and $\dot u_0 \in \Hu[\u_0]$ and $\alpha_0 \in L^2(\frictionBoundary)$.
\end{problem}
It is natural to start the evolution out of an equilibrium configuration, i.e.,
with an initial displacement $\u_0$  that solves the stationary problem
\begin{equation} \label{eq:INITU}
u_0 \in \Hu[\u_0]:\qquad  b(\u_0, \v) = \ell(\v) \qquad \forall \v \in \Hu[\u_0].
\end{equation}
In our numerical experiments to be reported below,  \eqref{eq:INITU} 
is solved iteratively by a fixed point iteration over the geometric nonlinearity,
i.e., starting with $\u_0^0 = 0$, a new iterate $\u_0^{\nu+1}$ is computed as the
(up to tangential rigid body motions) unique solution of  the corresponding linear problem on $\Hu[\u_0^\nu]$.

To our knowledge, existence and uniqueness of solutions of Problem~\ref{eq:WEAK}  is widely open.
In case of unilateral frictional contact with a rigid foundation and  Dieterich's law \eqref{eq:DIETERICH},
long-time existence of solutions was established by Pipping~\cite{Pipping17}.

%
\section{Semi-discretization in time}
Utilizing Rothe's method \cite{kavcur1986,rothe_1930}, we first perform a time discretization of Problem~\ref{eq:WEAK}
leading to a sequence of continuous spatial problems to be (approximately) solved in each time step.
To this end, the time interval $[0,T_0]$ is partitioned into time steps $0=t_0< \cdots < t_N=T_0$ 
with given step size $\tau_n= t_{n+1}-t_n>0$, and we write $\tau = \tau_n$ for notational convenience.

\subsection{Rate problem with given state}
We first consider the rate problem \eqref{eq:WEAKRATE} for given state $\alpha \in L^2(\Gamma^F)$.
Following \cite{PippingSanderKornhuber15},  we apply the classical Newmark scheme
\begin{equation} \label{eq:NEWMARK}
\begin{array}{rl}
	\dot{\u}_n & = \dot{\u}_{n-1} + \tfrac{\tau}{2} \left( \ddot{\u}_{n-1} + \ddot{\u}_{n} \right) \\
	\u_n & = \u_{n-1} + \tau \dot{\u}_{n-1} + \left(\tfrac{\tau}{2}\right)^2 \left( \ddot{\u}_{n-1} + \ddot{\u}_{n} \right)
	\end{array},
	\qquad n=1,\dots,N,
\end{equation}
which is well-known to be energy-conserving, consistent with second order, and unconditionally stable~\cite{Hairer-2006}.
Utilizing \eqref{eq:NEWMARK}, we eliminate
 \begin{equation} \label{eq:ELIMINATE}
\begin{array}{rl}
	\ddot{\u}_n & = \tfrac{2}{\tau} \left(\dot{\u}_n - \dot{\u}_{n-1}\right) -  \ddot{\u}_{n-1}, \\
	u_n       & = u_{n-1 }+ \frac{\tau}{2}(\dot u_{n} + \dot u_{n-1}), 
	\end{array},
	\qquad n = 1,\dots,N,
\end{equation}
from \eqref{eq:WEAKRATE} at fixed time $t=t_n$ and freeze the solution dependence in the closed-fault condition
and  in the friction law  at $\u_{n-1}$ to obtain the spatial variational inequality
\begin{equation} \label{eq:RATESPATIAL}
  \dot u_n \in \Huold: \qquad   a_n(\dot \u_n, \v-\dot{\u}_n) +  \Phi^{\u_{n-1}}(\v, \alpha) - \Phi^{\u_{n-1}}( \dot \u_n , \alpha)     \geq \ell_n(\v-\dot{\u}_n), 
  \qquad \forall \v \in \Huold,
\end{equation}
for $n=1, \dots, N$. Here,  we have set
\begin{equation} \label{eq:ACONT}
    a_n (\v, \w) = \tfrac{2}{\tau} (\rho \v, \w) +  a (\v, \w) +  \tfrac{\tau}{2} b(\v, \w) 
\end{equation}
with $(\cdot,\cdot)$ denoting the canonical scalar product in $L ^2(\Omega)$ and 
\begin{equation*}
    \ell_n(\v) = \ell(\v) + (\rho \ddot{\u}_{n-1}, \v) + \tfrac{2}{\tau} (\rho \dot{\u}_{n-1}, \v) 
    -  \tfrac{\tau}{2} b( \dot{\u}_{n-1},\v)  -  b(\u_{n-1},\v).
    \end{equation*}
Note that $\ddot u_0$  is not given as an initial condition in the continuous Problem~\ref{eq:WEAK}.
Assuming initial acceleration  towards equilibrium,  $\ddot u_0$ is therefore computed from the auxiliar problem
\begin{equation} \label{eq:INITACC}
\ddot \u_0 \in H_0: \qquad (\rho \ddot \u_0,  \v) + b(\u_0, v)  = \ell(\v)  \qquad \forall \v \in H_0.
\end{equation}

Note that the jump terms $\penu{\cdot}$, the contact
boundary $\frictionBoundaryU$, and the contact mapping $\pi^\u$
are all taken with respect to the last deformed state
$(\Id + \u_{n-1})(\frictionBoundary)$ of the contact boundary.
This eliminates the geometric nonlinearity associated with
large (relative) deformations of the contact boundary,
and we are left with the variational inequality \eqref{eq:RATESPATIAL}
on the  affine subspace $\Huold$ of $H_0$ to be solved in each time step..

As $a_n (\cdot, \cdot)$ is symmetric and positive definite and $\Phi(\cdot,\alpha)$ is convex, the variational inequality \eqref{eq:RATESPATIAL} 
can be equivalently written as the minimization problem
\begin{equation} \label{eq:RATEMIN}
  \dot u_n \in \Huold: \qquad   \cJ ( \dot \u, \alpha) \leq     \cJ (\v, \alpha) \qquad \forall  \v \in \Huold
\end{equation} 
for the corresponding energy functional
\begin{equation*}
\cJ (\v, \alpha) = \tfrac{1}{2} a_n (\v,\v) + \Phi^{\u_{n-1}} (\v, \alpha) - \ell_n (\v).
\end{equation*}
The following lemma~\cite[Theorem 6.49]{fonseca2007modern} will be useful to show existence and uniqueness of a solution.
\begin{lemma}\label{lem:LEMMA1}
Assume that $g:  \frictionBoundary \times \R^s \to \R$, $s \in \N$, is a non-negative function, 
such that $g(x,\cdot)$ is lower semicontinuous for almost all $x \in \frictionBoundary$.
Then the induced functional
\[
\int_{\frictionBoundary} g(x, \cdot) \; dx : \quad L^2(\frictionBoundary) \to \R \cup \{+ \infty\}
\]
is lower semicontinuous.
\end{lemma}

The convex functional $\Phi^{\u_{n-1}}(\cdot, \alpha)$ defined in \eqref{eq:RATEFUN} for the Dieterich--Ruina model  \eqref{eq:DRF}
is proper and lower semicontinuous by Lemma~\ref{lem:LEMMA1}.
Furthermore, by the assumptions on $\Ten{A}$ and $\Ten{B}$,
the bilinear form $(\Ten{A}(x)+ \frac{\tau}{2}\Ten{B}(x))( \cdot) \colon (\cdot)$ on the symmetric $\dim \times \dim$ matrices
is symmetric and uniformly elliptic with respect to $x\in \Omega$.
The following existence result therefore follows from Korn's second inequality and~\cite[Lemma~4.1]{Glowinski1984}.
%
\begin{proposition} \label{prop:RATESPATIALSOL}
Let $f \in L^2(\Omega)$ and $f^N \in L^2(\neumannBoundary)$.
Assume that $u_{n-1}$, $n = 1,\dots,N$,  avoids self-penetration  so that the contact mapping  $\pi^{u_{n-1}}$ and thus $\Huold$ are well-defined.
Then  the spatial rate problem \eqref{eq:RATESPATIAL} has  a unique solution  ${\dot u}_n\in \Huold$
and any given state $\alpha \in L^2(\frictionBoundary)$.
\end{proposition}
As a consequence of Proposition~\ref{prop:RATESPATIALSOL}, 
the solution operator  $R: L^2(\frictionBoundary) \to \Huold$,
\begin{equation} \label{eq:RATEOPERATOR}
  L^2(\frictionBoundary) \ni \alpha \mapsto R(\alpha) = \dot u_n\in \Huold,
\end{equation} 
of the spatial rate problem \eqref{eq:RATESPATIAL} is well-defined,
if no self-penetration occurs in preceding time steps. 
This is a strong assumption, as the contact conditions are taken explicitly.

\subsection{State problem with given rate}
Discretizing the state problem \eqref{eq:WEAKSTATE} with given velocity $\dot \u_n \in \Huold$
by the implicit Euler method and freezing the state law at $\u_{n-1}$
yields the variational inequality 
\begin{align} 
\begin{split}
	\alpha_n \in L^2(\frictionBoundary): \quad 	
  \left( \alpha_{n} , \beta - \alpha_{n} \right)_{L^2( \frictionBoundary)}  
  + \tau \Psi^{\u_{n-1}} & (\beta, \dot \u)  - \tau  \Psi^{\u_{n-1}} (\alpha_n, \dot \u )\\
	&  \geq  \left( \alpha_{n-1} , \beta - \alpha_n\right)_{L^2 (\frictionBoundary)}
\quad \forall \beta \in L^2 (\frictionBoundary)
\end{split}  \label{eq:STATESPATIAL}
\end{align}
which can be equivalently expressed as  the minimization problem
\begin{equation}\label{eq:STATEMIN}
\alpha_n \in L^2(\frictionBoundary): \quad     \cE(\alpha_n, \dot \u_n)  \leq   \cE(\beta, \dot \u_n) \qquad \forall \beta \in L^2 (\frictionBoundary)
\end{equation}
for the associated energy
\[
\cE(\beta, \dot \u) =  \tfrac{1}{2} (\beta, \beta)_{L^2 (\frictionBoundary)} + \tau \Psi^{\u_{n-1}} (\beta, \dot \u) - (\alpha_{n-1}, \beta)_{L^2 (\frictionBoundary)}.
\]
Both for  Dieterich's law \eqref{eq:DIETERICHLAW} and  Ruina's law  \eqref{eq:RUINALAW},
the functional $\Psi(\cdot, \dot u)$ defined in \eqref{eq:STATEFUN} 
is convex, proper and, by Lemma~\ref{lem:LEMMA1}, lower semicontinuous 
(see the proof of \cite[Proposition~4.4]{PippingSanderKornhuber15} for details).
Hence, existence and uniqueness again follows from~\cite[Lemma~4.1]{Glowinski1984} .
\begin{proposition} \label{prop:STATESPATIALSOL}
Both for Dieterich's law \eqref{eq:DIETERICHLAW} and  Ruina's law  \eqref{eq:RUINALAW}, the spatial state problem \eqref{eq:STATESPATIAL} has  a unique solution $\alpha_n\in L^2(\frictionBoundary)$ for $n = 1,\dots,N$  and any given velocity $\dot u \in \Huold$.
\end{proposition}
Proposition~\ref{prop:STATESPATIALSOL} gives rise to the solution operator $S: \Huold \to  L^2(\frictionBoundary)$,
\begin{equation} \label{eq:STATEOPERATOR}
  \Huold \ni \dot u  \mapsto S(\dot u) = \alpha_n \in L^2(\frictionBoundary),
\end{equation} 
of the spatial state problem  \eqref{eq:STATESPATIAL}.

\subsection{Coupled spatial problem}
Combining \eqref{eq:RATESPATIAL} and \eqref{eq:STATESPATIAL}, the time discretization of Problem~\ref{eq:WEAK} now reads as follows.
\begin{problem}[Semi-discretization in time] \label{prob:COUPLEDSPATIAL}
    Find $\dot \u_n \in \Huold$ and $\alpha_n \in L^2 (\frictionBoundary)$ satisfying
    \begin{align*}
    a_n(\dot \u_n, \v-\dot{\u}_n) +  \Phi^{\u_{n-1}}(\v, \alpha_n) - \Phi^{\u_{n-1}}( \dot \u_n , \alpha_n)    & \geq \ell_ n(\v-\dot{\u}_n), 
      &&  \forall \v \in \Huold\\
      \left( \alpha_{n} , \beta - \alpha_{n} \right)_{L^2( \frictionBoundary)} + \tau \Psi^{\u_{n-1}} (\beta, \dot \u_n) - \tau  \Psi^{\u_{n-1}} (\alpha_n, \dot \u_n )
		&\geq  \left( \alpha_{n-1} , \beta - \alpha_n\right)_{L^2 (\frictionBoundary)}
      &&  \forall  \beta \in L^2 (\frictionBoundary)
    \end{align*}
 for $n=1, \dots,N$ with 
 $\ddot u_n$  computed from \eqref{eq:ELIMINATE} and the auxiliary problem \eqref{eq:INITACC},
and  given initial conditions $u_0, \dot u_0 \in \Hu[\u_0]$, $\alpha_0 \in L^2 (\frictionBoundary)$.
\end{problem}


%

Recall, that $\Huold$ is well-defined only if $u_{n-1}$ avoids self-penetration, because the contact map $\pi^{u_{n-1}}$
is not available otherwise. This drawback could be overcome by  introducing an approximate contact map $\tilde{\pi}$
as in the spatial discretization below.

For an unilateral version of Problem~\ref{prob:COUPLEDSPATIAL}, 
i.e., a subduction zone with rigid foundation and Dieterich's law, such difficulties do not occur and
existence and uniqueness have been shown in~\cite[Proposition 3.6.]{Pipping2014}
based on Banach's fixed point theorem.
In case of Ruina's law, existence (but possibly no uniqueness)
was established utilizing Schauders fixed point theorem 
(see \cite[Corollary 3.8.]{Pipping2014} or \cite[Theorem 5.14]{PippingSanderKornhuber15}).


\section{Discretization in time and space}
For each $i=1, \dots, \bodyCount$
we assume that the subdomain $\Omega_i$ is polygonal and denote by $\cT_i$ a triangulation, i.e., a shape-regular, simplicial partition, 
of $\Omega_i$
with vertices $\cN^*_i$. 
We introduce the associated vector-valued, linear finite element space
\begin{equation*}
    \cS_{i} = \left\lbrace \v \in C(\Omega_i)^{\dim} \, | \, \v \text{ is linear on  all\;} T \in \cT_i  \text{ and } \v|_{\dirichletBoundary[i]}=0\right\rbrace .
\end{equation*}

Assuming that the Dirichlet boundary $\dirichletBoundary[i]$ is resolved by
$\cT_i$, 
we  define the set of nodes $\cN_i = \cN^*_i  \Setminus \overline{\dirichletBoundary[i]}$.
The resulting partition $\cT = \bigcup_{i=1}^{\bodyCount}\cT_i$ of $\Omega$  leads to the associated product space
\[
\cS = \cS_1\times \cdots \times \cS_{\bodyCount} = \op{span}\{\lambda_p \e_j \;|\; p\in \cN, \; j=1, \dots,\dim \} \subset H_0
\]
with the nodes $\cN =  \bigcup_{i=1}^{\bodyCount}\cN_i$, the nodal basis functions  $\lambda_p\in \cS_i$, $p \in \cN_i$, 
and the unit vectors $\e_j \in \Rd$, $j=1,\dots,\dim$. We emphasize that  the triangulations $\cT_i$ and $\cT_{i+1}$
do not need to match at the common interface  $\frictionBoundary[i,i+1]$, $i=1,\dots\bodyCount-1$,
in the sense that the sets $\cN_i \cap \frictionBoundary[i,i+1]$ and $\cN_{i+1} \cap \frictionBoundary[i,i+1]$ in general do not coincide.
For ease of notation, we even assume $\cN_i \cap \cN_j = \emptyset$ for $i \neq j$ so that 
we do not need to distinguish shared nodes and unshared nodes in the following presentation.

\subsection{Mortar discretization of rate problem with given state}
We consider the spatial rate problem \eqref{eq:RATESPATIAL} with given state $\alpha \in L^2(\frictionBoundary)$.
In order to incorporate  non-penetration and tangential friction along the fault system $\frictionBoundary$, 
we first introduce its triangulation 
\[
\cT^F= \bigcup_{i=1}^{\bodyCount-1} \cT_i^F, \qquad  \cT_i^F=\{F = T \cap \frictionBoundary[i,i+1]\;|\; T \in \cT_i\},
\]
with the nodes  $\cN^F= \bigcup_{i=1}^{\bodyCount-1} \cN_i^F$, $\cN_i^{F}= \cN_i \cap \frictionBoundary[i,i+1]$,
together with the corresponding trace space
\[
\cS^F = (\v_1, \dots,\v_{\bodyCount - 1}) \subset L^2(\frictionBoundary)^{\dim}, \qquad \v_i \in \cS_i^{F} = \cS_i|_{\frictionBoundary[i,i+1]},
\] 
spanned by the nodal basis $\lambda_p \e_j|_{\frictionBoundary}$, $p\in \cN^F$, $j=1, \dots,\dim$.
Note that the triangulation $\cT_i^F$  and the associated finite element trace space $\cS_i^F$ on $\frictionBoundary[i,i+1]$ 
are inherited  from $\Omega_i$ (the bottom non-mortar side), 
and do not coincide with corresponding traces from $\Omega_{i+1}$ (the top mortar side).

Mimicking the continuous case, the discretization of  non-penetration condition  and  friction  law
is based  on an  approximation  $\piDisc: \BOT{\frictionBoundaryDisc} \to \TOP{\frictionBoundaryDisc}$
of the contact mapping of $\pi^{\u_{n-1}}$ from the preceding time step with  corresponding approximations 
$\BOT{\frictionBoundaryDisc}\subset \frictionBoundary$, and $\TOP{\frictionBoundaryDisc}\subset \frictionBoundary$
of $\BOT{\frictionBoundaryU[\u_{n-1}]}$, and $\TOP{\frictionBoundaryU[\u_{n-1}]}$, respectively.
The approximations $\piDisc$ and $\BOT{\frictionBoundaryDisc}$  come into play, because
the top and bottom interfaces of the deformed subdomains $(\Id + \u_{n-1,i})(\Omega_i)$ 
may not match due to discretization errors that arise from enforcing non-penetration for $\u_{n-1}$.
In the following, we assume that the non-mortar contact
boundary $\BOT{\frictionBoundaryDisc}$ is resolved by a subset of the fault triangulation $\cT^F$.
We refer to~\cite{BastianBuseSander2010,Sander2011:psurface} and the references cited therein for
algorithms to compute such approximations of non-matching discrete intersections.

In analogy to \eqref{eq:CONTJUMP}  the jump of  $\v \in \cS$ across the discrete deformed contact boundary  is then defined by
\begin{align*}
  \penDisc{\v}
    =
    \BOT{\v} - \TOP{\v} \circ \tilde{\pi} &\qquad \text{on }\BOT{\frictionBoundaryDisc}.
\end{align*}
In the spirit of \cite{wohlmuth2003monotone}, the non-penetration condition appearing in \eqref{eq:SOLSPACE} 
and the (tangential) jumps appearing  in the functionals $\Phi^{\u_{n-1}}$ and $\Psi^{\u_{n-1}}$ of the Dieterich--Ruina model
will be incorporated in a weak sense with respect to a discrete test space
spanned by dual mortar basis functions
as introduced by Wohlmuth~\cite{wohlmuth2000mortar}.

To this end,
we first introduce the set of non-mortar contact nodes
\begin{align*}
  \NodesNM = \cN^F \cap \overline{\BOT{\frictionBoundaryDisc}}
\end{align*}
as well as the deformed contact set
\begin{align*}
  \CDisc = (\Id + \u_{n-1,B})(\BOT{\frictionBoundaryDisc}),
\end{align*} 
and we work with the pullback
of the $L^2(\CDisc)$ inner product  to $\BOT{\frictionBoundaryDisc}$
\begin{align*}
  \innerCDisc{v,w} = (v \circ (\Id + \u_{n-1,B})^{-1}, w \circ (\Id + \u_{n-1,B})^{-1})_{L^2(\CDisc)}.
\end{align*}
Then the dual mortar basis functions $\varphi_q$, $q \in \NodesNM$,  are defined to be piecewise linear on $\cT^F$,
have the same support as $\lambda_q|_{\frictionBoundary}$, and  satisfy the bi-orthogonality property
\begin{align*}
  \innerCDisc{\lambda_p|_{\frictionBoundary}, \varphi_q} = \delta_{p,q}
  \qquad \forall p,q \in \NodesNM \qquad (\text{Kronecker-}\delta).
\end{align*}
Note that dual mortar functions are typically discontinuous and therefore not contained in $\cS^F$.
We refer to~\cite{wohlmuth2000mortar,wohlmuth2011variationally} for details about the construction.

We now define the linear projection $\Pi: \cS \to \cS$, componentwise according to
\begin{align*}
  (\Pi v)_j  = \Pi v_j = v_j - \sum_{p \in \NodesNM}\innerCDisc{\penDisc{v_j}, \varphi_p}\lambda_p, \qquad j=1,\dots d.
\end{align*}
Observe that here and in the following we denote by $\Pi$ the projection of both scalar and $\R^d$-valued functions.
The projection $\Pi$ gives rise to the direct splitting
\begin{align*}
  \cS = \cV \oplus \cW
\end{align*}
of $\cS$ into the image $\cV = \op{im}\Pi$  and the kernel $\cW = \op{ker}\Pi$ of $\Pi$.
Utilizing the fact that $\penDisc{\lambda_p} = \lambda_p$ for $p \in \NodesNM$,
we find that
these spaces can be written as
\begin{align*}
  \cV
  &= \left\{ \v \in \cS  \midst \innerCDisc{\penDisc{\v},\varphi_p}=0 \;\forall p\in \NodesNM \text{ and }j=1, \dots, d  \right\}\\
  &= \op{span}\left\{\mu_p \e_j \midst  \mu_p=  \Pi \lambda_p, \; p \in \cN\Setminus \NodesNM \text{ and } j=1,\dots,d \right\}
\end{align*}
and
\begin{align*}
  \cW
  &= \left\{\v \in \cS \midst \v(p)=0 \; \forall p \in \cN \Setminus \NodesNM \right\}
  = \op{span}\left\{\lambda_p \e_j \midst  p \in \NodesNM \text{ and } j=1,\dots,d \right\}.
\end{align*}
Thus, $\cV$ consists of all functions in $\cS$ which are weakly continuous across
the deformed contact boundary with respect to the
pullback $L^2$ scalar product   $\innerCDisc{\cdot,\cdot}$ and the associated dual mortar space
(notice that point-wise continuity does not hold  in general).

Correspondingly, the  basis functions  $\mu_p = \Pi \lambda_p$ spanning $\cV$
are  the usual hat functions $\lambda_p$ on the mortar side,
which are extended in a weakly continuous way to the non-mortar side,
while the basis functions $\lambda_p\e_j$ of $\cW$ involve the usual hat
functions  $\lambda_p$, $p \in \NodesNM$ on the non-mortar side
which drop down to zero across $\BOT{\frictionBoundaryDisc}$.

Both normal and tangential jumps of $v \in \cV$ are weakly zero.
As a consequence, both (weak) normal and tangential  jumps can be represented in terms of the incremental space $\cW$.
To this end, we define a nodal approximation $\normal_{\cS}$ of the normal to the deformed contact set $\CDisc$ by
\begin{equation} \label{eq;NODALNORMAL} 
  \normal_{\cS} = \sum_{p \in \NodesNM} \normal_p \lambda_p|_{\BOT{\frictionBoundaryDisc}},  \qquad \normal_{p} = \frac{\sum_{F \in \cT^F_p} \normal_F}{| \sum_{F \in \cT^F_p} \normal_F  |},
\end{equation}
where $\cT^F_p$ denotes the set of simplices $F \in \cT_i^F$ with common vertex $p\in \NodesNM$, 
and $\normal_F$ is an approximate normal to the deformed face
$\BOT{(\Id + \u_{n-1})}(F)$, e.g. the average of the normal on $\BOT{(\Id + \u_{n-1})}(F)$.
We also introduce the  approximate tangent space
$T_p \CDisc =  (\op{span} \{\normal_p \})^\perp \subset \R^d$
to the deformed contact set $\CDisc$ associated with the nodal approximate normal $\normal_p$ in $p\in \NodesNM$.
Similar to the continuous case, the discrete normal field $\normal_{\cS}$ to $\CDisc$
is parametrized over $\BOT{\frictionBoundaryDisc}$.

Next, we further split $\cW$  into its normal and tangential part
\begin{align*}
  \cW = \cW_n \oplus \cW_t,
\end{align*}
with
\begin{align*} 
  \cW_n & = \left\{ \lambda_p x \midst p \in \NodesNM, x \in \op{span}\{\normal_p\} \right\}, &
  \cW_t & = \left\{ \lambda_p x \midst p \in \NodesNM, x \in T_p \CDisc \right\}.
\end{align*}
Excluding normal jumps, we now define a (non-conforming) finite element counterpart of the solution space $\Huold$ according to
\begin{equation} \label{eq:DISCSOLSPACE}
  \Suold = \cV \oplus \cW_t.
\end{equation}
Such a mortar approach to non-penetration has been suggested and first analyzed in~\cite{wohlmuth2000mortar}.

The splitting suggests the unique decomposition
\begin{equation} \label{eq:SEPREP}
\v = \v_{\cV} + \v_{\cW}, \qquad \v \in \Suold
\end{equation}
 where 
\[
  \v_{\cV}=\Pi v = \sum_{\cN\Setminus \NodesNM} v(p) \mu_p \in \cV
\]
 and
\[
  \v_{\cW} = (\Id-\Pi)\v = \sum_{p \in \NodesNM} \pen{\v}_p \lambda_p \in \cW_t, 
\]
denoting the weak nodal jump of $v$ at $p \in \NodesNM$ by
\begin{align*}
  \pen{\v}_p = (\v - \Pi(\v))(p) = (\innerCDisc{\penDisc{\v_j},\varphi_p})_{j=1}^d \in \R^d.
\end{align*} 
The nodal vectors $\pen{\v}_p$ clearly satisfy $\pen{\v}_p \cdot \normal_p=0$ for all $p \in \NodesNM$.
Hence, $v_{\cW}$ can be regarded as  a nodal approximation of the  tangential jump of $v$ along $\CDisc$ pulled back to 
$\BOT{\frictionBoundaryDisc}$.
Inserting this approximation into \eqref{eq:RATEFUN} and 
replacing the integrand $\phi(\penu{\v}, \alpha)$  by its nodal interpolation in $\cS^F$,
we obtain the approximate functional $\Phi_{\cS}: \Suold \to \R$, 
\begin{equation} \label{eq:NODPHI}
  \Phi_{\cS}(\v, \alpha) =  \sum_{p \in \NodesNM} \phi_p(\pen{\v}_p), \qquad 
  \phi_p(\pen{\v}_p) = \phi(\pen{v}_p,\alpha(p))\int_{\BOT{\frictionBoundaryDisc}}\lambda_p\; ds.
\end{equation}
At this point, $\alpha$ is required to be continuous in a neighborhood of each node $p \in \NodesNM$ to guarantee  that nodal interpolation makes sense.

The mortar discretization of the rate problem \eqref{eq:RATESPATIAL} 
with given, sufficiently regular state $\alpha\in L^2(\frictionBoundary)$  now reads as follows
\begin{equation} \label{eq:DISCRATESPATIAL}
 \dot u_{n,\cS} \in \Suold: \qquad   a_n(\dot \u_{n,\cS}, \v-\dot{\u}_{n,\cS}) +  \Phi_{\cS}(\v, \alpha) - \Phi_{\cS}( \dot \u_{n,\cS} , \alpha)     \geq \ell_ {n,\cS}(\v-\dot{\u}_{n,\cS}), 
 \qquad \forall \v \in \Suold,
\end{equation}
for $n=1, \dots, N$. Here,  the bilinear form $a_n(\cdot,\cdot)$ is taken from \eqref{eq:ACONT} and we have set
\begin{equation*}
    \ell_{n,\cS}(\v) = \ell(\v) + (\rho \ddot{\u}_{n-1,\cS}, \v) + \tfrac{2}{\tau} (\rho \dot{\u}_{n-1,\cS}, \v) 
    -  \tfrac{\tau}{2} b( \dot{\u}_{n-1,\cS},\v)  -  b(\u_{n-1,\cS},\v)
    \end{equation*}
with $\ddot \u_{n-1,\cS}$, $\dot \u_{n-1,\cS}$, $\u_{n-1,\cS}$ taken from preceding time steps, by discrete analogues of \eqref{eq:ELIMINATE}, 
by suitable finite element approximations  $\dot \u_{0,\cS}$, $\u_{0,\cS}\in \cS$ of the initial conditions $\dot \u_0$, $\u_0 \in \Hu[\u_0]$,
or a finite element approximation $\ddot \u_{0,\cS}\in \cS$ of the auxiliary problem \eqref{eq:INITACC}.
Existence and uniqueness of discrete  spatial solutions $ \dot u_{n,\cS} \in \Suold$, $n=1,\dots,N$, follows 
under the same conditions and by the same arguments as in Proposition~\ref{prop:RATESPATIALSOL}.

This mortar approach  to (frictional)  non-penetration directly extends elastic frictional contact problems, i.e., to fault opening.
We refer to~\cite{krause2009nonsmooth,krause2002dirichlet} for further information and to~\cite{wohlmuth2011variationally} for a detailed survey.

\subsection{Piecewise constant discretization of state problem with given rate}
We consider the state problem \eqref{eq:STATESPATIAL} 
with a given deformation rate.
Let $\cC^F= \{C_p \subset \frictionBoundary\;|\;  p \in \cN^F\}$
be a dual partition of the triangulation $\cT^F$ of $\frictionBoundary$.
We introduce the subspace $\cB^F \subset L^2(\frictionBoundary)$ of functions that are  constant on each cell $C_p \in \cC^F$, $p\in \cN^F$,
and the resulting piecewise constant discretization
\begin{equation}  \label{eq:DISCSTATESPATIAL}
\begin{array}{rcl}
	\alpha_{n,\cB} \in \cB^F: \quad 	
  \left( \alpha_{n,\cB}  , \beta - \alpha_{n,\cB}  \right)_{L^2( \frictionBoundary)} +
   \tau \Psi _{\cB}(\beta, \dot \u) - \tau  \Psi_{\cB}(\alpha_{n,\cB}, \dot \u)
		\geq & \\
		  \left( \alpha_{n-1,\cB} , \beta - \alpha_{n,\cB}\right)_{L^2 (\frictionBoundary)}
 \;\;&  \; \forall \beta \in \cB^F 
\end{array}
\end{equation}
of the state problem \eqref{eq:STATESPATIAL} with given $\dot \u \in \Suold$. 
Here, the nodal approximation  $ \Psi_{\cB}: \cB^F \to \R$,
\begin{equation} \label{eq:DISCSTATEFUN}
  \Psi_{\cB}( \beta,  \dot u) =   \sum_{p \in \NodesNM} \psi (\beta(p), |\pen{\dot \u}_p|) \abs{C_p}, \qquad \beta \in \cB^F,
\end{equation}
of the functional $ \Psi(\cdot, \dot u)$ is obtained in the same way as 
the nodal approximation $\Phi_{\cS}(\cdot, \alpha)$ of $\Phi(\cdot, \alpha)$  in \eqref{eq:NODPHI}.

Existence and uniqueness of discrete  spatial solutions $ \alpha_{n,\cB} \in \cB^F$, $n=1,\dots,N$, follows 
in the same way as in Proposition~\ref{prop:STATESPATIALSOL}.

\subsection{Fully discretized coupled spatial problem}
Combining \eqref{eq:DISCRATESPATIAL} and \eqref{eq:DISCSTATESPATIAL}, the  discretization of 
the coupled  Problem~\ref{eq:WEAK} in time and space  now reads as follows.

\begin{problem}[Discretization in time and space] \label{prob:DISCCOUPLEDSPATIAL}
    Find $\dot \u_{n, \cS} \in \Suold$ and $\alpha_{n,\cB} \in \cB^F$ satisfying
    \begin{align*}
a_n(\dot \u_{n,\cS}, \v-\dot{\u}_{n,\cS}) + 
  \Phi_{\cS}(\v, \alpha_{n,\cB} ) - \Phi_{\cS}( \dot \u_{n,\cS} , \alpha_{n,\cB} )   &  \geq \ell_ {n,\cS}(\v-\dot{\u}_{n,\cS})
 \qquad \forall \v \in \Suold\\
  \left( \alpha_{n,\cB}  , \beta - \alpha_{n,\cB}  \right)_{L^2( \frictionBoundary)} + \tau \Psi _{\cB}(\beta, \dot \u_{n,\cS}) - \tau  \Psi_{\cB}(\alpha_{n,\cB}, \dot \u_{n,\cS} )
		 & \geq  
		  \left( \alpha_{n-1,\cB} , \beta - \alpha_{n,\cB}\right)_{L^2 (\frictionBoundary)} \quad \forall \beta \in \cB^F   
    \end{align*}
 for $n=1, \dots,N$ with given initial conditions $u_{0,\cS}, \dot u_{0,\cS} \in \Suold$, $\alpha_{0,\cB} \in \cB^F$.
\end{problem}


Iterative solution of Problem~\ref{prob:DISCCOUPLEDSPATIAL}
can be obtained from the fixed point iteration
\begin{equation} \label{eq:DISCFIXEDP}
   \dot{\u}_{n,\cS}^{\nu+1} = R_{\cS} \left( \omega \alpha_{n,\cB}^{\nu+1} + (1- \omega) \alpha_{n,\cB}^{\nu} \right), \quad
    \alpha_{n,\cB}^{\nu+1} = S_{\cB} \left( \dot{\u}_{n,\cS}^{\nu} \right), \quad   \nu = 0,1, \ldots, 
\end{equation}
with initial iterate  $ (\dot{\u}_{n,\cS}^0 , \alpha_{n,\cB}^0) =  (\dot{\u}_{n-1,\cS} , \alpha_{n-1,\cB})$ 
and suitable relaxation parameter $\omega \in (0,1]$. 
Here, $S_{{\cB}}:  \Suold \to \cB^F$  and  $R_{{\cS}}:  \cB^F \to \Suold$  denote the solution operators of
the state problem with given rate \eqref{eq:DISCSTATESPATIAL} and 
the rate problem with given state \eqref{eq:DISCRATESPATIAL}, respectively.
Note that state functions $\alpha \in\cB^F$ are continuous in a neighborhood of  each node $p \in \cN^F$ and thus
satisfy the regularity assumptions made for nodal interpolation \eqref{eq:NODPHI}.

Extension of the convergence proof given in \cite{Pipping2014,PippingSanderKornhuber15}
for an unilateral version of Problem~\ref{prob:DISCCOUPLEDSPATIAL}
to the actual discretized multi-body problem is a subject of future research.


%
\section{Algebraic solution}

\subsection{Fixed point iteration and state problem with given rate}
The iterative solution of the  coupled Problem~\ref{prob:DISCCOUPLEDSPATIAL}  
is performed by the fixed point iteration \eqref{eq:DISCFIXEDP}.
The state problem with given rate~\eqref{eq:DISCSTATESPATIAL} arising in each iteration step
fully decouples into scalar algebraic problems for $\alpha_{n,\cB}(p)$ and  $p \in \cN^F$ 
that can be solved explicitly or, e.g., by bisection.

The rate problem with given state~\eqref{eq:DISCRATESPATIAL}, however, is a  discretized frictional contact problem,
and its iterative solution is more involved.

\subsection{Truncated non-smooth Newton multigrid for the rate problem with given state} \label{subsec:TNNMG}
We now concentrate on the robust and efficient algebraic solution of the
mortar-discretized rate problem  \eqref{eq:DISCRATESPATIAL} with given state $\alpha \in  \cB^F$.
First, recall that the splitting  \eqref{eq:SEPREP}  provides the basis representation
\begin{equation} \label{eq:HIERBAS}
  \v = \sum_{p \in \cN\Setminus \NodesNM} v_p\mu_p + \sum_{p \in \NodesNM} v_p \lambda_p 
\end{equation}
of all $v \in \Suold$ with coefficients $\v_p=\v(p)$, $p\in \cN\Setminus \NodesNM$ and  $\v_p=\pen{\v}_p$, $p\in\NodesNM$.
We identify each $\v \in \Suold$  with its coefficient vector $(\v_p)_{p\in \cN}$. 
Then the discrete nonlinear functional 
\[
  \Phi_{\cS}(\v, \alpha) = \sum_{p \in \NodesNM} \phi_p (v_p)
\]
introduced in \eqref{eq:NODPHI} has a separable structure in the sense that the  coefficients $v_p \in \R^d$ 
are decoupled with respect to the nonlinearity $\phi_p$.

The variational inequality  \eqref{eq:DISCRATESPATIAL} can be equivalently
rewritten as the minimization problem
\begin{equation} \label{eq:DISCRATESPATIALMIN}
\dot u_{n,\cS}\in \Suold:\qquad \cJ_\cS(\dot u_{n,\cS}) \leq \cJ_\cS(\v)\qquad \forall \v \in \Suold
\end{equation}
denoting
\[
  \cJ_\cS(\v)  = {\textstyle \frac{1}{2} }a_n(\v, \v) - \ell_{n, \cS}(\v) +  \sum_{p \in \NodesNM} \phi_p(v_p).
\]
This formulation  allows to construct and analyze globally convergent  nonlinear  
Gau\ss--Seidel relaxation methods~\cite{Glowinski1984}. Based on the splitting
\begin{align*}
  \Suold = \sum_{p \in \cN} \cV_p, \qquad
 \cV_p=
\begin{cases}
  \{ \mu_p x\; | \; x \in \R^d\}  &\text{for } p \in \cN \Setminus \NodesNM,  \\
  \{\lambda_p x\;|\; x \in T_p \CDisc \}   &\text{for } p \in \NodesNM
  \end{cases}
\end{align*}
and some enumeration $\cN = \{p_1, \dots, p_M\}$,
a new iterate is computed by successive subspace minimization:
Given an iterate $u$ set $w_0 = u$ and compute $w_i$, $i=1,\dots,M$
by solving 
\begin{equation} \label{eq:GSRELAX}
  w_i \in w_{i-1}+\cV_{p_i}: \quad  \cJ_\cS(w_i ) \leq \cJ_\cS(w) \quad \forall   w \in w_{i-1}+\cV_{p_i}, \qquad i=1, \dots, M,
\end{equation}
to obtain the new iterate $\overline{\u}=w_M$.
However, such iterative schemes are well-known to 
suffer from rapidly deteriorating convergence rates for decreasing mesh size.

The basic idea of \emph{Truncated Nonsmooth Newton Multigrid} (TNNMG) 
methods~\cite{Graeser2011,GraeserSackSander2009,GraeserSander2019}
is to complement nonlinear Gau\ss--Seidel smoothing \eqref{eq:GSRELAX}
by additional line search into the  Newton-type search direction $\newtonCorr$,
as obtained from the linear system
\begin{equation}\label{eq:LINEARIZATION}
\cJ_{\cS}''(\bar{u})|_{W(\bar{\u})\times W(\bar{u})} \newtonCorr = - \cJ_{\cS}'(\bar{u})|_{W(\bar{\u})}
\end{equation}
on a suitable subspace $W(\bar{\u})\subset \Suold$.
Accounting for non-smoothness of $\phi_p$, $p \in \NodesNM$, we select the reduced subspace
 \[
   W(\bar{\u}) = \cV + \op{span}\{ \lambda_p x \;|\; x \in T_p \CDisc \; |\pen{\bar{\u}}_p|\neq V_m(\alpha(p)), \; p \in  \NodesNM \}.
 \]
By freezing $\bar{u}(p)$ at those $p$ where $\phi_p(\bar{u}_p)$ is not smooth enough, the
restriction $\cJ_{\cS}|_{W(\bar{\u})\times W(\bar{u})} $  to $W(\bar{\u})$ is twice differentiable.
Note that global convergence of nonlinear Gau\ss--Seidel smoothing \eqref{eq:GSRELAX} is preserved by any correction $\rho \newtonCorr$
such that $\rho \in [0,\infty)$ is providing non-increasing energy 
\begin{equation} \label{eq:NONDEGENER}
\cJ_\cS (\bar{u} + \rho \newtonCorr)\leq \cJ_\cS (\bar{u}).
\end{equation}

In TNNMG methods, all  three substeps, i.e., nonlinear Gau\ss--Seidel relaxation \eqref{eq:GSRELAX}, 
evaluation of the Newton-type search direction \eqref{eq:LINEARIZATION}, 
and monotone line search \eqref{eq:NONDEGENER}, are typically performed inexactly.
A TNNMG  iteration step applied to a given iterate $\u^\nu \in \Suold$ thus reads as follows
\begin{equation} \label{eq:TNNMG}
 \bar{\u}_\nu = P (\u^{\nu}),
\quad 
  \newtonCorr^\nu = MG(\cJ_{\cS}''(\bar{u}_\nu)|_{W(\bar{\u}_\nu)\times W(\bar{u}_\nu)},  \cJ_{\cS}'(\bar{u}_\nu)|_{W(\bar{\u}_\nu)}),
 \quad 
 u^{\nu+1}= \bar{\u}_{\nu} + \rho (\bar{\u}_{\nu}, \newtonCorr^{\nu}) \newtonCorr^\nu, 
\end{equation}
with corresponding inexact solution operators  $P$, $MG$, and $\rho$.  

More precisely,  inexact Gau\ss--Seidel relaxation $P$  is obtained as follows.
For each convex $(d-1)$-dimensional non-smooth minimization problem \eqref{eq:GSRELAX}
on the local tangent space $\cV_{p_i}$ associated with the node $p_i \in \NodesNM$,
the quadratic part (i.e., the $i$-th diagonal block of the stiffness matrix 
corresponding to the bilinear form $a_n(\cdot,\cdot)$ and the basis representation \eqref{eq:HIERBAS})
is replaced by a scalar upper bound, e.g. its maximal eigenvalue.
Then the resulting problem is rotationally symmetric in $\cV_{p_i}$ and thus reduces to a
scalar problem that can be solved by bisection or even by an explicit formula (cf.~\cite[Example~5.2]{GraeserSander2019}).
The same local preconditioning approach is used for the linear $d$-dimensional problems
in $\cV_{p_i}$ for $p_i \in \cN \Setminus \NodesNM$.
We emphasize that global convergence of nonlinear Gau\ss--Seidel relaxation \eqref{eq:GSRELAX} 
is preserved in this way~\cite{GraeserSander2019}.

Since $\cJ_\cS$ is strongly convex, the coefficient matrix $\cJ_{\cS}''(\bar{u}_\nu)|_{W(\bar{\u}_\nu)\times W(\bar{u}_\nu)}$
of the linear problem \eqref{eq:LINEARIZATION} is symmetric and positive definite
on the subspace $W(\bar{\u}_\nu)$. Hence,  its inexact solution $MG$
can be simply performed by one or more steps of a standard linear multigrid method
with minor modifications to deal with the special basis used in~\eqref{eq:HIERBAS} and the restriction to $W(\bar{\u}_\nu)$.
For details on the choice of suitable coarse grid spaces or, equivalently,  suitable restriction and
prolongation operators we refer, e.g.,  to~\cite{GraeserSackSander2009,sander:multidim_coupling_knee:2008}.

Inexact line search providing a damping factor $\rho_\nu \in [0,\infty)$ 
that guarantees monotonically decreasing energy \eqref{eq:NONDEGENER}
is finally performed by bisection.

The following convergence result is obtained as a special case of the abstract result~\cite[Corollary~4.5]{GraeserSander2019}
 by making use of~\cite[Theorem~5.6 and Lemma~5.8]{GraeserSander2019} 
 to incorporate the  inexact pre-smoothing $P$ explained above.

\begin{proposition}
  For any initial iterate $u^0 \in \Suold$ the sequence $u^\nu\in \Suold$, $\nu =1,\dots$, 
  generated by the TNNMG method \eqref{eq:TNNMG}
  converges to the unique solution of the mortar-discretized rate problem  \eqref{eq:DISCRATESPATIAL} 
  with given state $\alpha \in  \cB^F$.
\end{proposition}

The same convergence result applies, if more than one nonlinear pre-smoothing step or
additional nonlinear post-smoothing is utilized.
Note that TNNMG methods allow for straightforward extensions to fault opening
by incorporating non-penetration into the nonlinear Gau\ss--Seidel
smoother and enforcing feasibility of coarse corrections $\newtonCorr^\nu$
by an additional projection step.

%
\section{Numerical experiments}

\subsection{General setup}
\subsubsection{Problem description}
\begin{figure}[ht]
  \centering
  \begin{minipage}[c]{0.45\textwidth}
    \includegraphics[width=\linewidth]{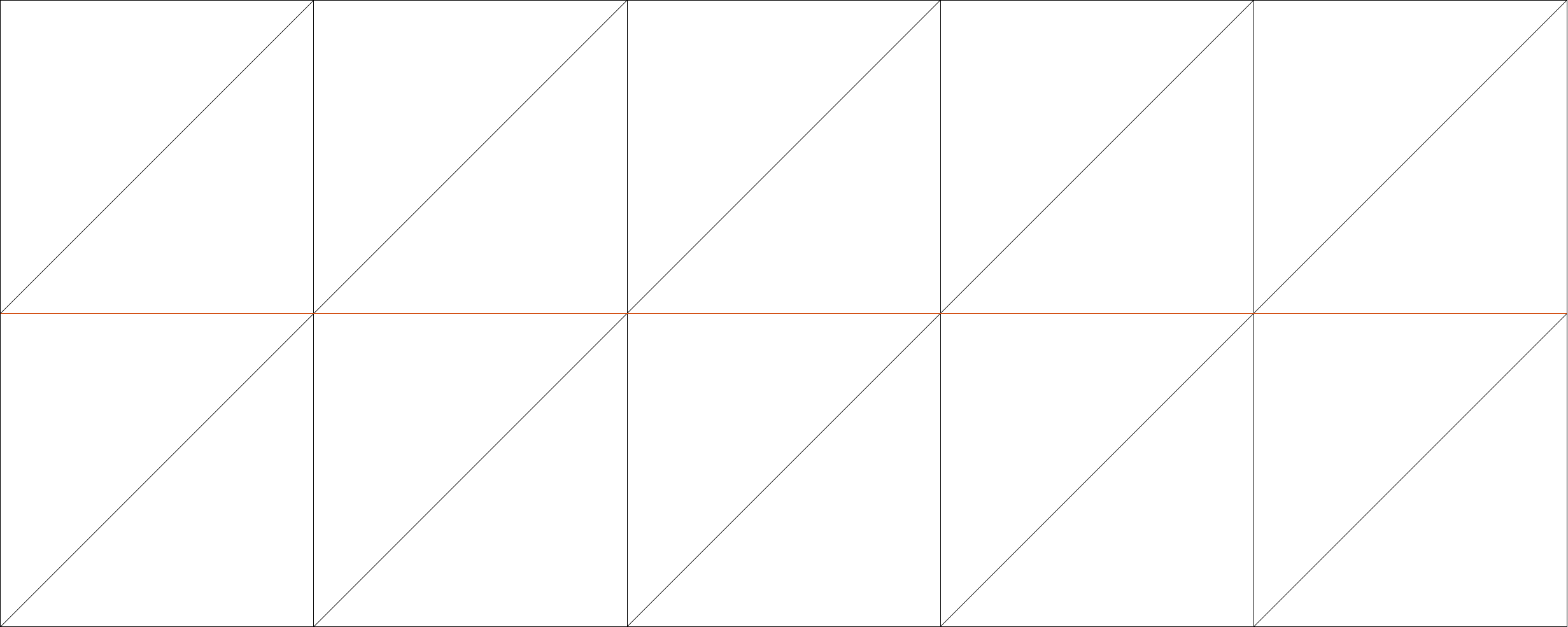}
  \end{minipage}
  \begin{minipage}[c]{0.45\textwidth}
    \includegraphics[width=\linewidth]{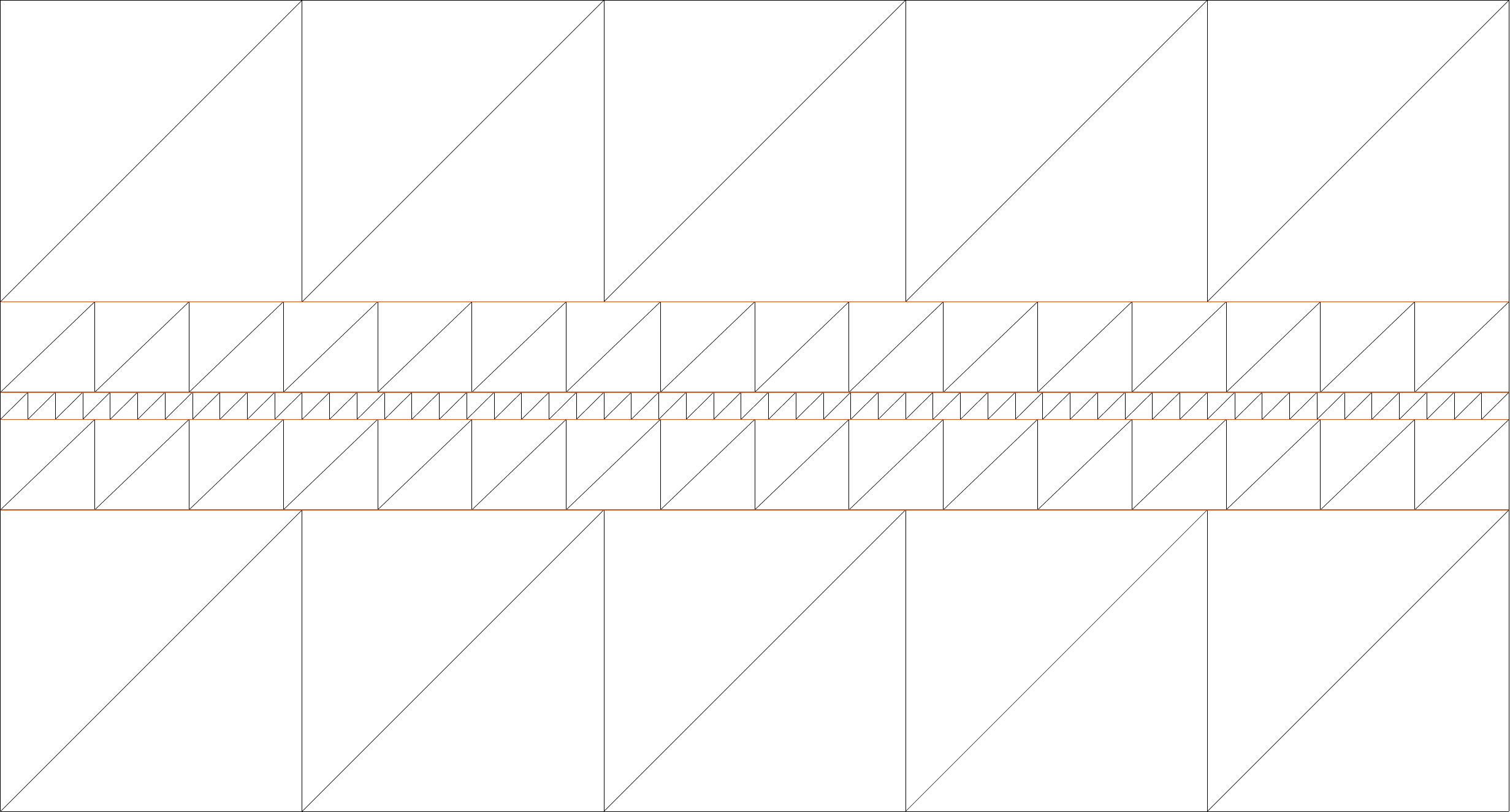}
  \end{minipage}
	\caption{Initial triangulations  ${\mathcal T}^{(0)}_{i}$, $i=1,\dots, \bodyCount$, for the spring slider with $\bodyCount=2$ bodies  (left)
	and a layered  fault system with $\bodyCount=5$ bodies (right).}  \label{fig:NUMGEO}
\end{figure}

%
In our two numerical experiments, we consider a rectangular deformable body in $d=2$ space dimensions 
that is decomposed into $\bodyCount=2$ (spring slider) or $\bodyCount=5$ (layered fault system) rectangular bodies 
by $1$ or $4$ planar faults, cf. Figure~\ref{fig:NUMGEO}.

In the spring slider experiment, the 2 bodies are both of the size $\SI{5}{\meter} \times \SI{1}{\meter}$.
They are associated with the reference domains $\Omega_1 = (-2.5, 2.5)\times (-1,0)$, $\Omega_2 = (-2.5, 2.5)\times (0,1)$
with the interface 
\[
\frictionBoundary = (-2.5, 2.5)\times \{0\}.
\] 
This setup corresponds to the one presented in \cite{PippingSanderKornhuber15}, but features a deformable instead of a rigid foundation. 

The bodies of the layered fault system have the size $\SI{5}{\meter} \times \SI{1}{\meter}$,  $\SI{5}{\meter}\times \SI{0.3}{\meter}$, $\SI{5}{\meter}\times \SI{0.09}{\meter}$,
$\SI{5}{\meter}\times \SI{0.3}{\meter}$, and $\SI{5}{\meter}\times \SI{1}{\meter}$.
They are associated with the reference domains
\begin{align*}
  \Omega_1 &= (-2.5, 2.5)\times (-1.345, -0.345),&
    \Omega_2 &= (-2.5, 2.5)\times (-0.345,-0.045),\\
  \Omega_3 &= (-2.5, 2.5)\times (-0.045, 0.045),&
    \Omega_4 &= (-2.5, 2.5)\times (0.045,0.345),\\
  \Omega_5 &= (-2.5, 2.5)\times (0.345, 1.345)
\end{align*}
with the interface  $\frictionBoundary$,
\[
\frictionBoundary = (-2.5, 2.5)\times \{-0.345\} \cup (-2.5, 2.5)\times \{-0.045\} \cup (-2.5, 2.5)\times \{0.045\} \cup (-2.5, 2.5)\times \{0.345\}.
\] 

The bodies consist of St.\ Venant Kirchhoff material and are subject to gravity, 
i.e., the body force is constant and  given by $f = - \rho g \cdot e_2$ with $g$ denoting the gravitational constant.
We impose homogeneous Neumann boundary conditions $f^N=0$
at the vertical  boundary   $\neumannBoundary$ of the associated reference configurations $\Omega_i$,  $i=1,\dots,\bodyCount$.
The system is fixed by homogeneous Dirichlet conditions $\u(\cdot,t)=\dot u(\cdot,t) = 0$, $0 \leq t \leq T_0$
at  the foundation $\dirichletBoundary[1]$.  At the upper Dirichlet boundary $\dirichletBoundary[\bodyCount]$, $I=2,\; 5$,
 the  condition $\dot{\u} (\cdot, t) = v_D \, \xi(t) \cdot e_1$ prescribes a  smooth transition 
 from zero  velocity to  constant loading speed $v_D = \SI[quotient-mode = fraction]{2e-4}{\meter / \second}$ with
\begin{equation} \label{eq:dirichlet_value}
    \xi (t) = 
    \begin{cases}
        \tfrac{1}{2} (1 - \cos(4 \pi t / T_0 )), & \text{if } t \leq T_0/10 \\
        1 & \text{otherwise}.
    \end{cases}
\end{equation}
At the  interfaces, $\frictionBoundary = \bigcup_{i=1}^{\bodyCount-1} \frictionBoundary[{i,i+1}]$ 
we impose rate-and-state friction conditions with Dieterich’s aging law.

The initial deformation $\u (\cdot, 0)$ is obtained by approximating the equilibrium configuration,
i.e., the solution of the stationary problem  \eqref{eq:INITU}, by one step of the associated fixed point iteration.
The initial velocity field is set to zero, which is consistent with the Dirichlet conditions,
and the initial state field is chosen to be $\alpha (\cdot, 0) = -10$ on $\frictionBoundary$. 

We consider the time interval $[0,T_0]$ with final time $T_0= \SI{60}{\second}$, and
the  remaining material parameters are given in Table \ref{tab:material_props}. 
	\begin{table}[ht]
        \centering
		\begin{tabular}{ll||ll}
        Bulk parameter & Value & Friction parameter & Value \\
        \hline
		Bulk modulus $E$ & $\SI{4.12e7}{\pascal}$ & ref. velocity $V_0$ & $\SI[quotient-mode = fraction]{1e-6}{\meter / \second}$ \\
		Poisson ratio $\nu$ & $\num{0.3} $ & ref. friction coeff. $\mu_0$ & $\num{0.6} $ \\
		mass density $\rho$ & $\SI[quotient-mode = fraction]{5e3}{\kilogram / \square\meter}$ & $a$ & $\num{0.010} $ \\
		gravity $g$ & $\SI[quotient-mode = fraction]{9.81}{\newton / \kilogram}$ & $b$ & $\num{0.015}$ \\
		& & charact. slip dist. $L$ & $\SI{1e-5}{\meter}$
		\end{tabular}
        \caption{Material parameters}
        \label{tab:material_props}
	\end{table}
	

\subsubsection{Discretization and algebraic solution}


\begin{figure}[ht]
  \setlength{\unitlength}{1cm}
  \centering
  \begin{minipage}[c]{0.48\textwidth}
    \includegraphics[width=\linewidth]{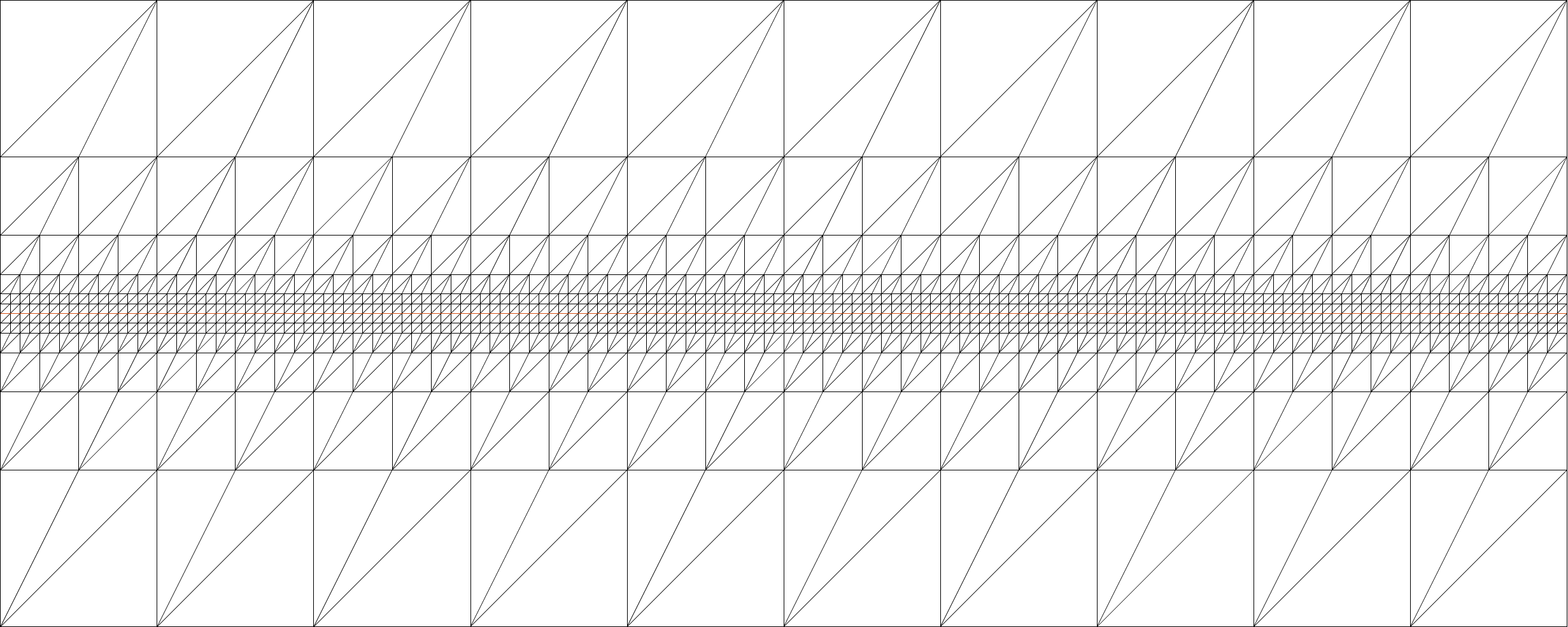}
  \end{minipage}
  \begin{minipage}[c]{0.48\textwidth}
    \includegraphics[width=\linewidth]{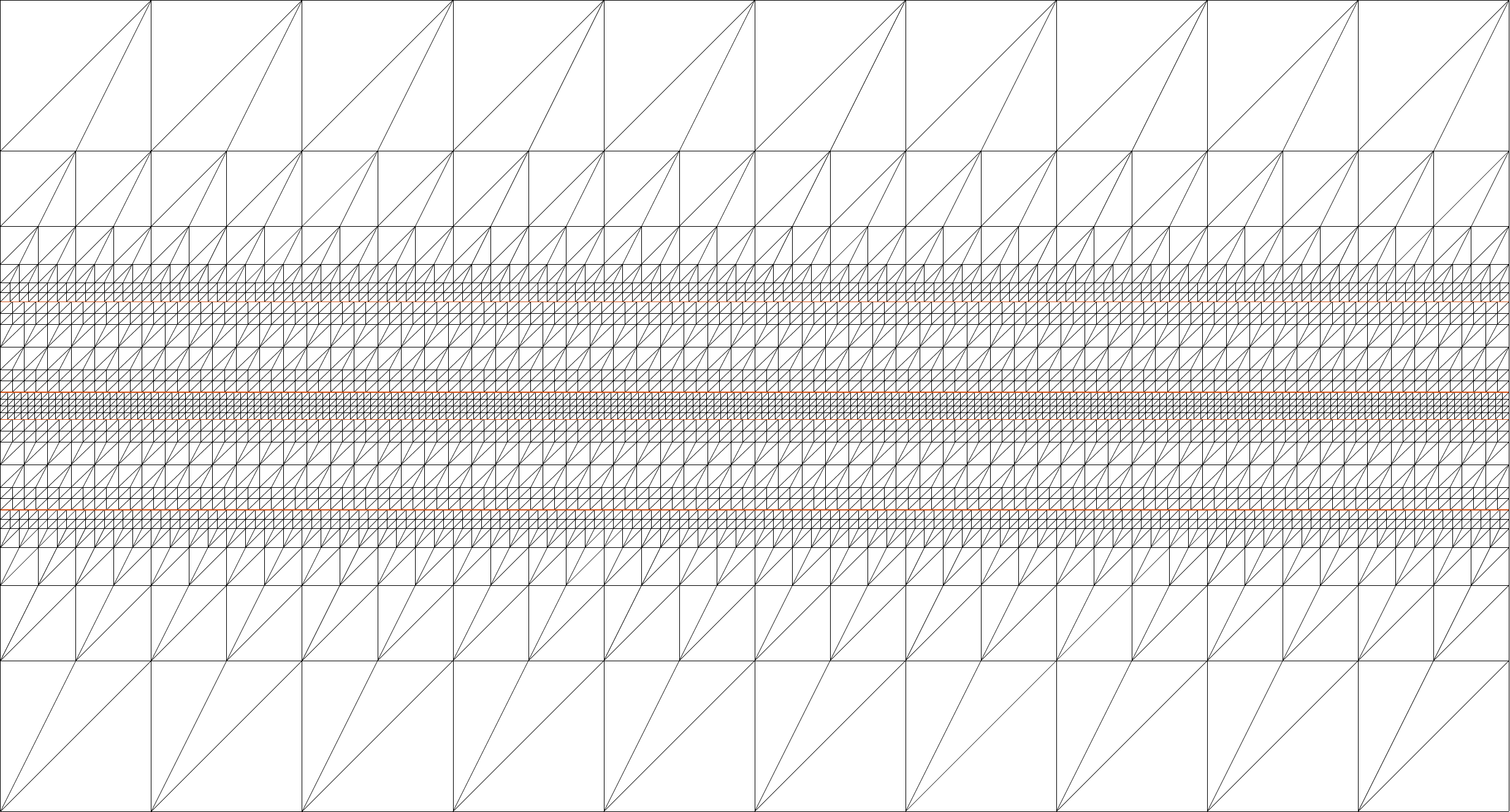}
  \end{minipage}
  \caption{Adaptively refined final triangulations  $T \in \cT^{(K)}_i$, $i=1,\dots, \bodyCount$ with $K=5$ and $1274$ vertices for the spring slider (left) and $K=5$ and $4057$ vertices for the layered fault system (right).}  \label{fig:NUMGEOFINAL}
\end{figure}

In order to efficiently resolve strongly varying dynamics ranging from slow interseismic loading to fast coseismic periods,
corresponding time step sizes are automatically selected according to  the following adaptive strategy.
For given approximate solution at $t_n$ of the coupled spatial Problem~\ref{prob:COUPLEDSPATIAL} at time $t_n \in [0,T_0)$,  
as computed by the old time step size $\tau_{n-1}$, 
we choose $\tau_{n}^* = \tau_{n-1}$ for $n\geq 1$ and $\tau_{-1}= 10^{-4} \, T_0$ as an initial guess for the new time step size $\tau_n$.
Then, we compute  approximate solutions 
$({\dot u}^{(1)}_{n+1}, {\alpha}^{(1)}_{n+1})$ at $t_{n} + 2 \tau_{n}^* $ by  one step with step size $2\tau_{n}^*$ 
and $({\dot u}^{(2)}_{n+1}, {\alpha}^{(2)}_{n+1})$ by two time steps with step size $\tau_{n}^*$.
If the criterion 
\begin{equation} \label{eq:TIMECRIT}
    \Vert   {\alpha}^{(1)}_{n+1} -  {\alpha}^{(2)}_{n+1}\Vert_{L^2 (\frictionBoundary)} \leq 
    \SI[parse-numbers=false]{\delta_{\tau}}{\meter^{1/2}}
\end{equation}
holds with a suitable threshold $\delta_{\tau}$, then we allow for coarsening:
With the new guess $\tau_{n}^*:= 2\tau_{n}^*$ 
the above procedure is repeated until \eqref{eq:TIMECRIT} is violated and we set $\tau_n := \tau_{n}^*/2$ in this case.
If the criterion  \eqref{eq:TIMECRIT} is already violated by the initial guess $\tau_{n}^* = \tau_{n-1}$,
then we require refinement: Successive bisection $\tau_{n}^*:= \tau_{n}^*/2$ is applied until  \eqref{eq:TIMECRIT} is met
and we set $\tau_n := \tau_{n}^*$ in this case.
The threshold  $\delta_{\tau}$ is selected in accordance with the accuracy of the inner fixed point iteration to be specified below.

The spatial problems occurring in each time step are discretized with respect to triangulations $\cT_i= \cT_i^{(K)}$,
as resulting from $K$ refinement steps  applied to initial triangulations $\cT^{(0)}_i$ of the subdomains $\Omega_i$, $i = 1, \dots, I$.
Note that the associated hierarchy of finite element spaces is utilized for the algebraic TNNMG solver to be specified later on.
For the spring slider and the layered fault system, the initial triangulations  $\cT^{(0)}_i$, $i=1,\dots, \bodyCount$, 
are shown in the left and in the right picture of Figure~\ref{fig:NUMGEO}, respectively.

For both geometries,  refinement is concentrated  at the interfaces  by the following adaptive procedure.
Starting with 
$\mathcal{T}_i^{(0)}$,
we perform regular (red) refinement of all  triangles $T\in  \cT^{(k)}_i$, $k\geq 0$, with diameter $h_T$  violating the criterion
\begin{equation} \label{eq:REFCRIT}
h_T < (1 + 80 \, d(T, \frictionBoundary) \, h_{\min}.
\end{equation} 
Here, $d(T, \frictionBoundary)$ stands for the distance of $T$ to $\frictionBoundary$
and $h_{\min} = \SI{6.25}{\centi\meter}$. Then, triangles with two or three bisected edges
as emerging through this procedure are also refined regularly until only triangles with no or with only one bisected edge are left
from ${\mathcal T}^{(k)}$.
The latter ones are then refined by connecting the midpoint of this edge with the opposite vertex to obtain  conforming refined triangulations
$\cT^{(k+1)}_i$.
In order to preserve shape regularity, these (green) closures are removed in advance of the next refinement step~\cite{bank1983some}.
Refinement terminates with $K=k$,
once the criterion \eqref{eq:REFCRIT} is met by all triangles $T \in \cT^{(k)}_i$ and all $i=1,\dots, \bodyCount$.

The resulting final triangulations are depicted in Figure~\ref{fig:NUMGEOFINAL}.
For the spring slider  (left), the final triangulations  $\mathcal{T}_i^{(K)}$ 
are resulting from  $K=5$ adaptive refinement steps,  have  $1274$ vertices in total, 
and 
$\SI[parse-numbers=false]{4.4}{\centi\meter}\leq h_T \leq \SI[parse-numbers=false]{70.8}{\centi\meter}$ 
holds for the diameters $h_T$ of all $T \in \mathcal{T}_i^{(K)}$, $i=1, 2$.
For the layered fault system (right),  the final triangulations  $\mathcal{T}_i^{(K)}$ 
are obtained after $K=5$ refinement steps,  have $4057$ vertices in total, 
and $\SI{3.2}{\centi\meter}\leq h_T \leq \SI{70.8}{\centi\meter}$ holds for the diameters 
$h_T$ of all $T \in \mathcal{T}_i^{(K)}$,  $i=1, \dots,5$.

 %
 In the fixed point iteration \eqref{eq:DISCFIXEDP},  providing the decoupling of rate and state,
 we use the relaxation parameter $\omega = 1/2$. The iteration is stopped, once the criterion 
\begin{equation} \label{eq:fpiStop}
    \Vert \alpha_{n,\cB}^{\nu} - \alpha_{n,\cB}^{\nu - 1} \Vert_{L^2 (\frictionBoundary)} 
    \leq 10^{-1} \SI[parse-numbers=false]{\delta_{\tau}}{\meter^{1/2}}
\end{equation}
is satisfied. Here, the parameter $\delta_{\tau}$ is the same as in the  time step selection criterion \eqref{eq:TIMECRIT}.
This choice aims at comparable accuracy of fixed point iteration and time stepping, 
and the actual value $\delta_{\tau}=10^{-5}$ is motivated by systematic trial and error, cf.~\cite[Subsection 3.3]{pipping2016efficient}.

The algebraic solution of the discrete state problem \ref{eq:STATESPATIAL} with Dieterich's aging law  and given rate 
is approximated by pointwise bisection. The iteration  is stopped, once the  error in each node is uniformly 
bounded by the threshold $\SI{e-12}{\meter^{1/2}}$.

Starting with the final iterate  ${\dot u}^0_{n,\cS} = {\dot u}^{\nu_{\rm stop}}_{n-1,\cS}$
from the preceding time step, the algebraic solution of the discrete rate problem~\eqref{eq:DISCRATESPATIAL} 
with given state is performed  by a Truncated Nonsmooth Newton Multigrid (TNNMG) method as described in Subsection \ref{subsec:TNNMG}. 
In each  iteration step, the (truncated) linear correction is obtained by $5$  steps of a 
classical multigrid V-cycle with $3$ pre- and $3$ post-smoothing steps.
Here, we utilize the grid hierarchy provided by successive refinement described above.
The iteration is terminated, once the stopping criterion
\begin{equation} \label{eq:MGSTOP}
\| {\dot u}^{\nu} - {\dot u}^{\nu-1}\|_{n} \leq \SI{e-8}{\watt^{1/2} \meter^{1/2}} 
\end{equation}
is satisfied with the time-dependent energy norm $\|\cdot\|_n = a_n(\cdot,\cdot)^{1/2}$ and  $a_n(\cdot,\cdot)$ defined in \eqref{eq:ACONT}.
This stopping criterion is selected to reduce the  error of the inner multigrid iteration some orders of magnitude below 
the error of the outer fixed point iteration which is intended to be in the range of the discretization error.

The discretization and algorithms are implemented using the Dune
framework~\cite{BastianEtAl2019_dune} making use of the dune-grid-glue library~\cite{BastianBuseSander2010}
for the mortar coupling.

\subsection{Spring slider}

\begin{figure}[ht]
	\centering
    \includegraphics[width=\linewidth]{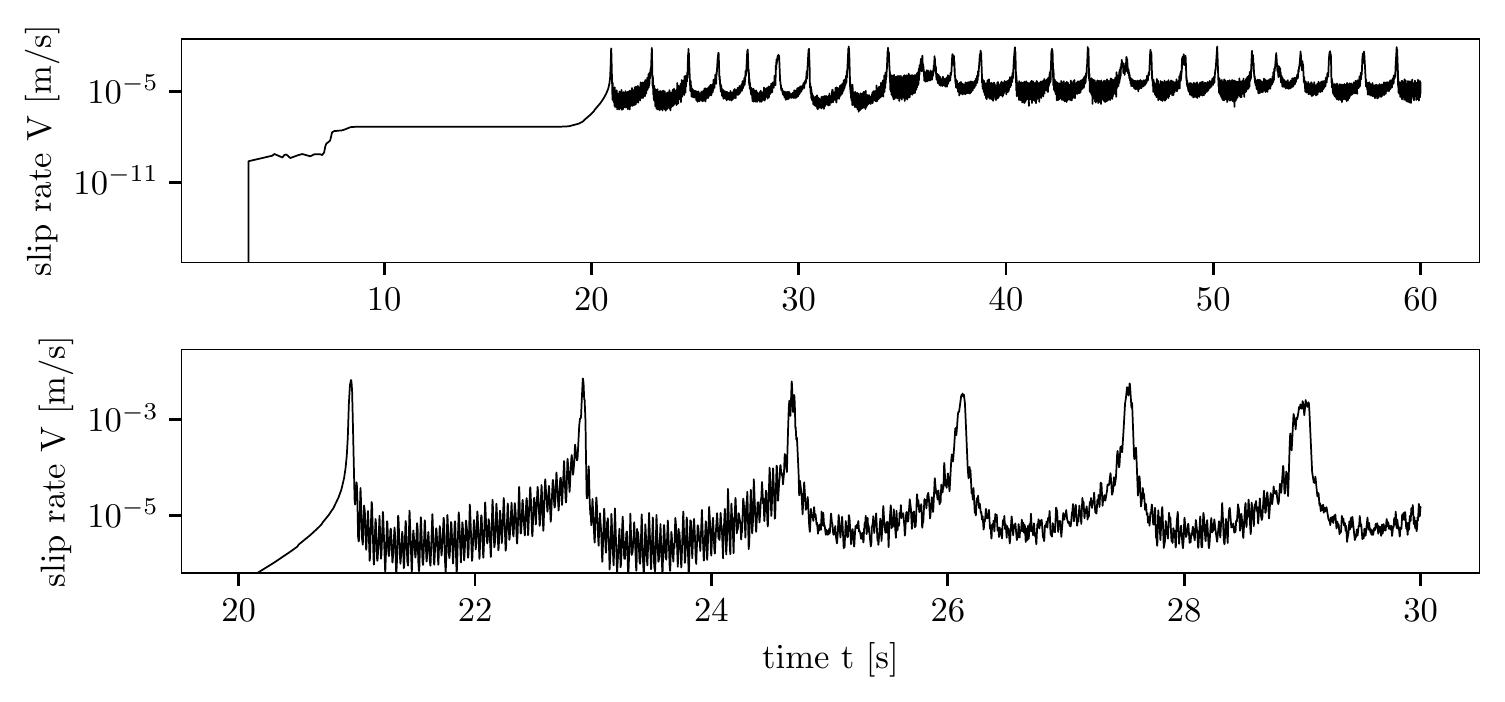}
    \vspace*{-25pt}
    \caption{Spring slider:
    Evolution of the mean value of relative velocities over $\frictionBoundary$  from
    the initial loading phase to slip events (top) and  zoom into the  slip events (bottom).}
    \label{fig:SPRISSLIPS}
\end{figure}

 
\begin{figure}
	\centering
    \includegraphics[width=\linewidth]{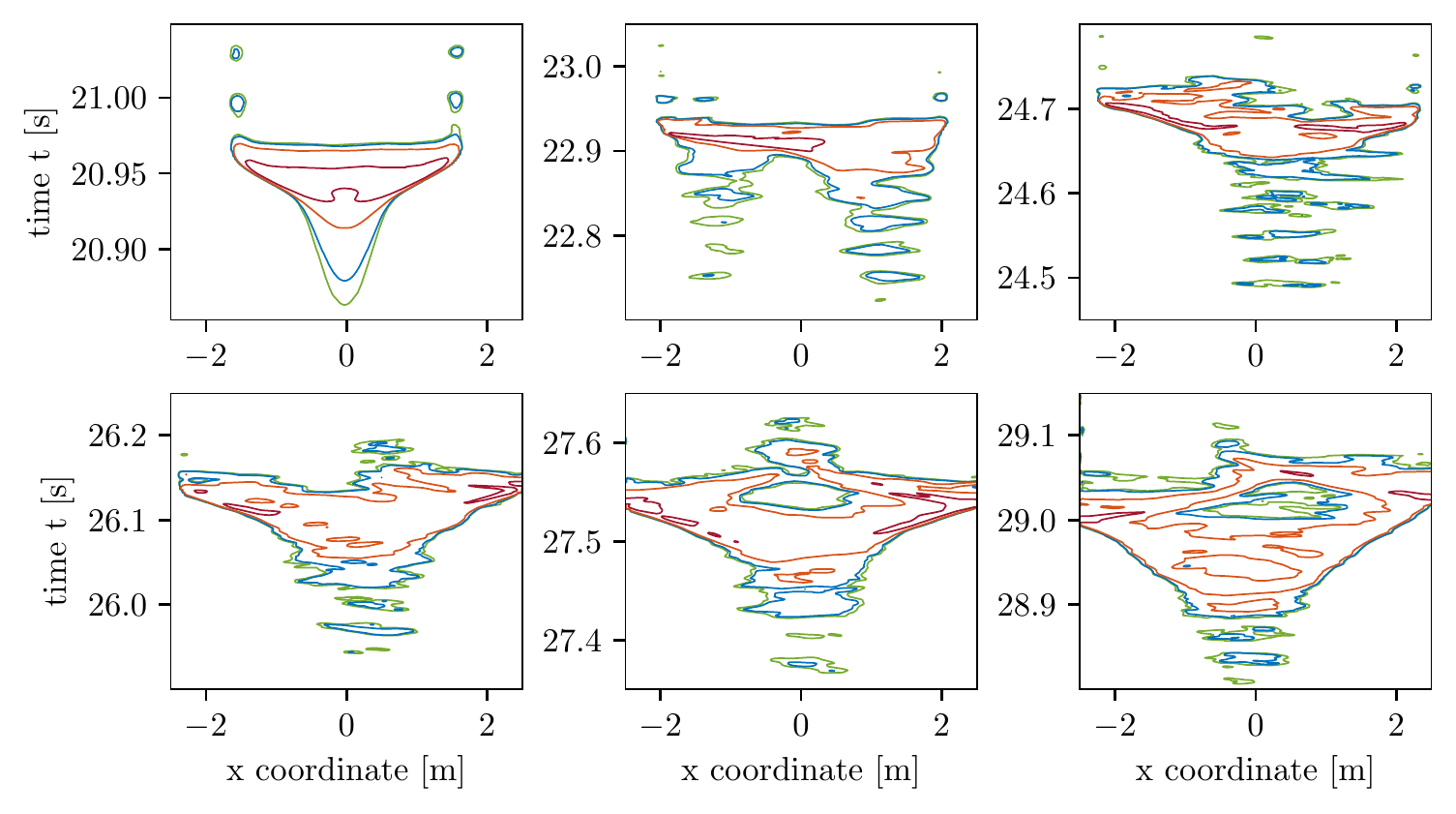}
    \vspace*{-25pt}
    \caption{Spring slider:
    Level lines ($\SI{e1}{\micro\meter / \second}$ \textcolor{plot-green}{green}, $\SI{e2}{\micro\meter / \second}$ \textcolor{plot-blue}{blue}, $\SI{e3}{\micro\meter / \second}$ \textcolor{plot-orange}{orange}, $\SI{e4}{\micro\meter / \second}$ \textcolor{plot-red}{red}) of the approximate relative velocity along  $\frictionBoundary$
    over  time intervals associated with the first 6 slip events}
    \label{fig:SPRISLEVEL}
\end{figure}

\subsubsection{Simulation results}
In order to illustrate the behavior of the deformed bodies along the fault  $\frictionBoundary= (-2.5,2.5)\times \{0\}$,
the top picture of Figure~\ref{fig:SPRISSLIPS}  shows the 
mean value of the approximate relative velocitiy $|\;[\dot u_{n,\cS}]^{u_{n-1,\cS}} |$ over $\frictionBoundary$
for the corresponding time instants $t_n$, $n= 1, \dots, N= 107659$. 
After a loading phase of about  $\SI{20}{\second}$, we observe 26 almost periodic peaks in the relative velocity,
indicating the occurrence of corresponding  slip events. 
Observe that periodicity is slightly perturbed in comparison with related numerical results for a rigid foundation~\cite{PippingSanderKornhuber15}.
A zoom into the first 6 slip events as shown in the bottom picture of Figure~\ref{fig:SPRISSLIPS}
reveals a highly oscillatory behavior of the approximate velocity. This is  partly due to the well-known 
lack of stability of the Newmark scheme~\cite{krause2012presentation},
but also occurs for  the highly dissipative backward Euler method~\cite{podlesny20} and sufficiently fine time steps.

Figure~\ref{fig:SPRISLEVEL} shows the level lines of  approximate relative velocities 
along the fault $\frictionBoundary$ (horizontal axis)
evolving over  time  (vertical axis) for six time intervals 
associated with the first six slip events (top left to bottom right).
The first (bilateral)  slip event originates from  the midpoint of  $\frictionBoundary$,
the next one has  two symmetric precursors on the left and right hand side,
and the third one features foreshocks with a slight emphasis towards  the right hand side of $\frictionBoundary$, 
but ruptures the entire fault nonetheless.
We then observe a sequence of  three further  bilateral  slip events beginning in the middle  of the fault.

All slip events are preceded by small foreshocks indicating the origin of the later event,
as well as small aftershocks  typically occurring at the right or left side of $\frictionBoundary$ or in its very center.
Note that these results considerably differ from related computations for a subduction zone
with rigid foundation that showed pure periodic behavior~\cite{pipping2016efficient}.

\subsubsection{Adaptive time stepping and  performance of the algebraic solver}
\begin{figure}
	\centering
    \includegraphics[width=\linewidth]{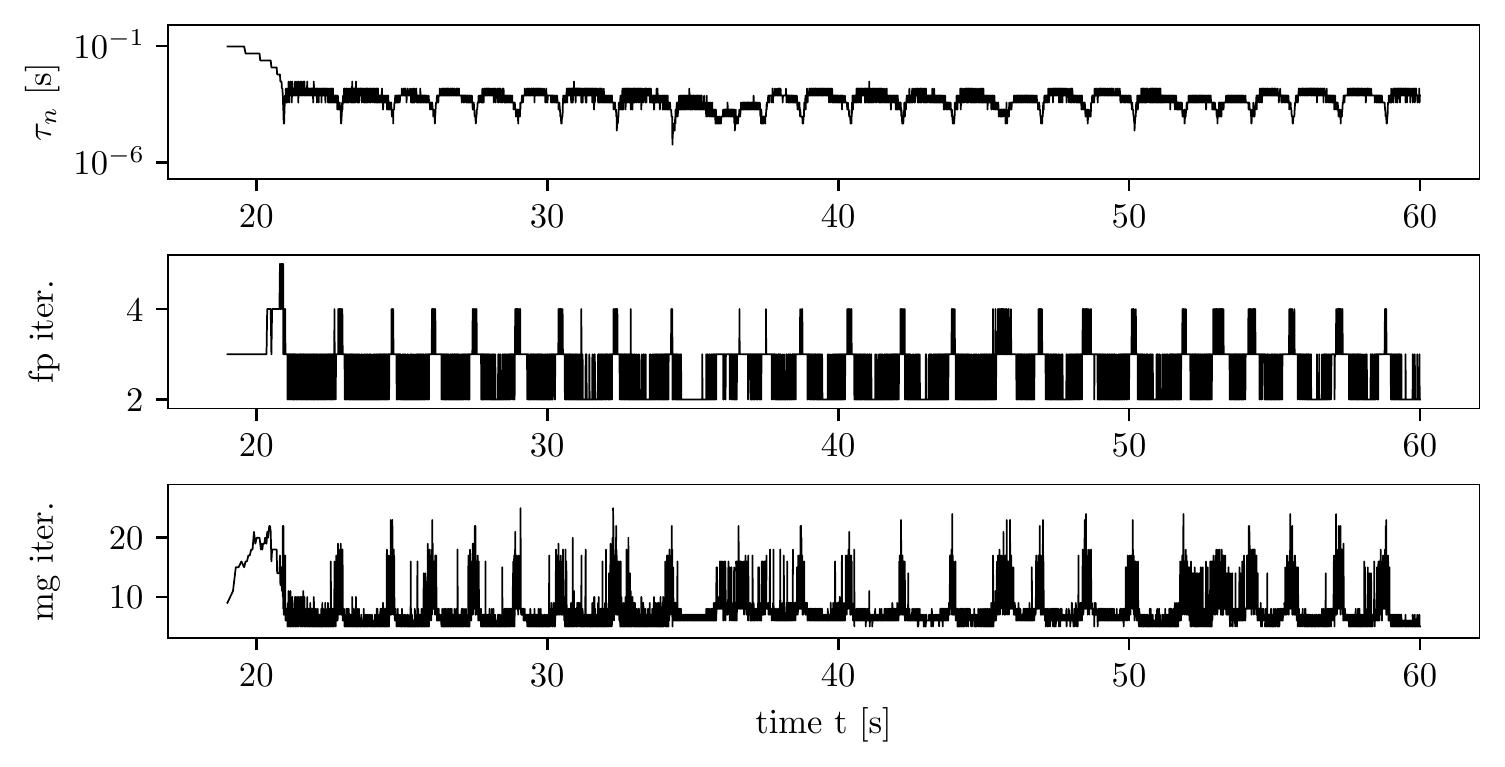}
    \vspace*{-25pt}
    \caption{Spring slider:
    Adaptive time step selection and  performance of the algebraic solver}
    \label{fig:SPRISSOLVER}
\end{figure}
We now describe the performance of adaptive time step selection and 
of the algebraic solver consisting of the fixed point iteration~\eqref{eq:DISCFIXEDP}
for the decoupling of rate and state and the TNNMG method for the rate problem as explained in Subsection~\ref{subsec:TNNMG}.
The upper picture of Figure~\ref{fig:SPRISSOLVER} shows the automatically selected time step sizes $\tau_n$ over  
corresponding time instants $t_n$ taken from the time interval that begins shortly before the end of the initial loading phase.
Observe that the occurrence of slip events is nicely reflected by the reduction of the time step size by about 2 orders of magnitude.
According to the second picture,
usually 2 - 4  fixed point iterations are required to match the stopping criterion~\eqref{eq:fpiStop}
for the actual spatial  problem with adaptively selected time step.
The third picture shows the sum of all inner multigrid iterations as  needed  to reach the stopping criterion \eqref{eq:MGSTOP}
in each of these outer fixed point iteration steps.
This sum often, but not always, increases and decreases with the number of required outer fixed point iterations  
and is ranging from about $5$ to  $29$.

\subsection{Layered fault system}
\begin{figure}
	\centering
    \includegraphics[width=\linewidth]{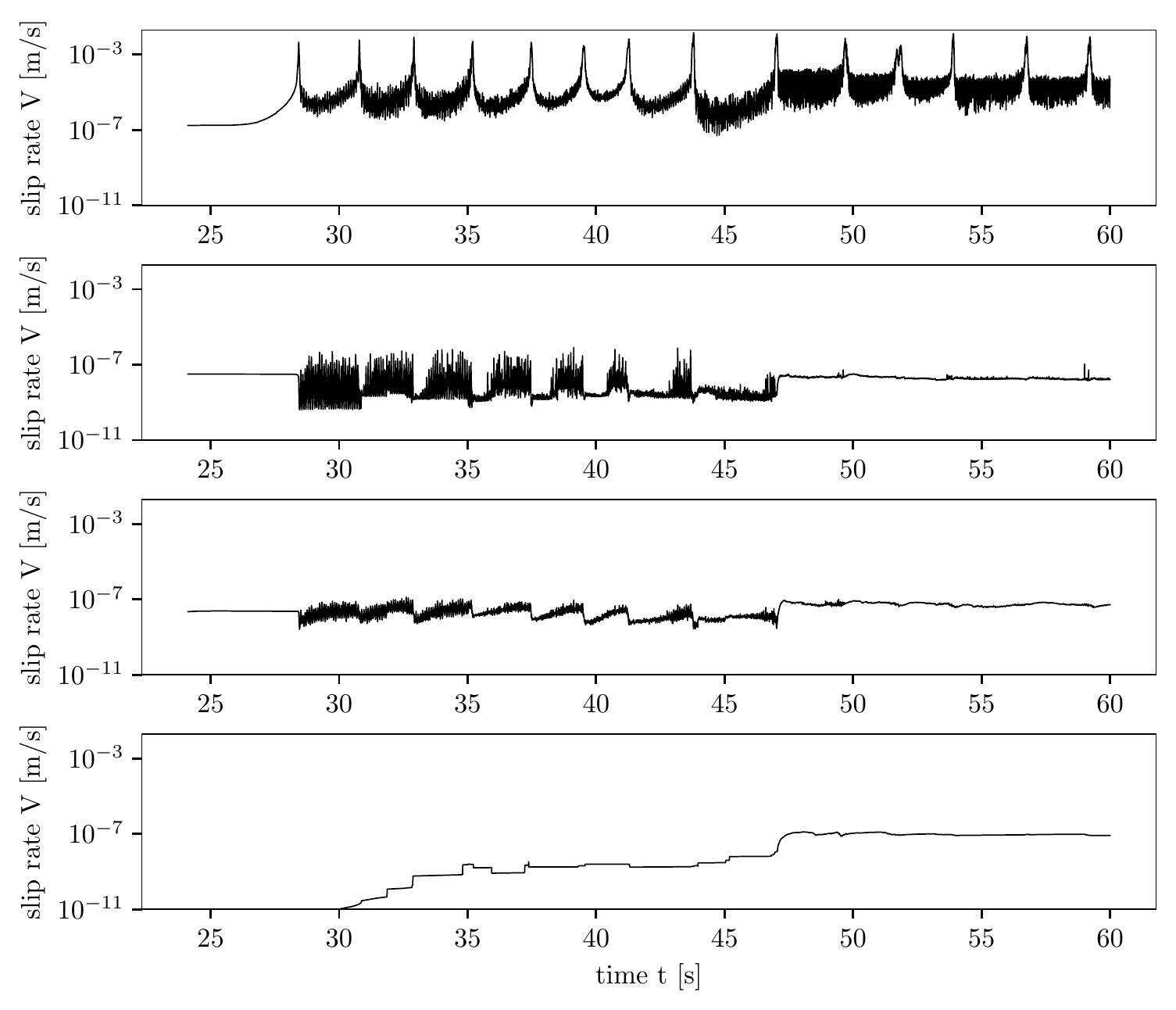}
    \vspace*{-25pt}
    \caption{Layered fault system:
    Evolution of the mean value of relative velocities over the faults $\frictionBoundary[i,i+1]$, $i=1,2,3,4$,
from $\frictionBoundary[4,5]$ (top) to $\frictionBoundary[1,2]$ (bottom).}
    \label{fig:MULLSLIPS}
\end{figure}

\begin{figure}
	\centering
    \includegraphics[width=\linewidth]{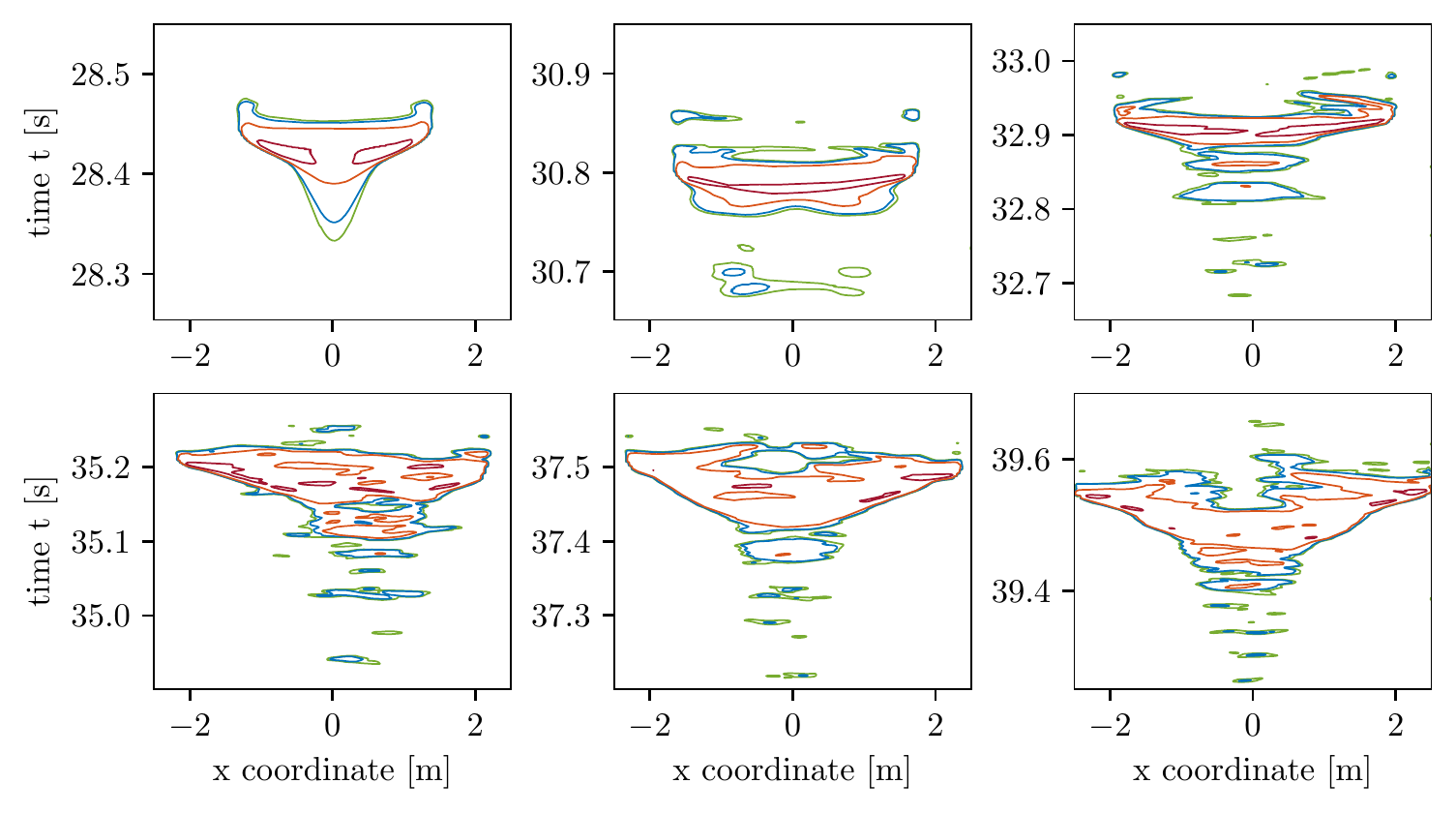}
    \vspace*{-25pt}
    \caption{Layered fault system:
    Level lines ($\SI{e1}{\micro\meter / \second}$ \textcolor{plot-green}{green}, $\SI{e2}{\micro\meter / \second}$ \textcolor{plot-blue}{blue}, $\SI{e3}{\micro\meter / \second}$ \textcolor{plot-orange}{orange}, $\SI{e4}{\micro\meter / \second}$ \textcolor{plot-red}{red}) of the approximate relative velocity along $\frictionBoundary[34]$ over  time intervals associated with the first 6 slip events.}
    \label{fig:MULLLEVEL}
\end{figure}

\begin{figure}
	\centering
    \includegraphics[width=\linewidth]{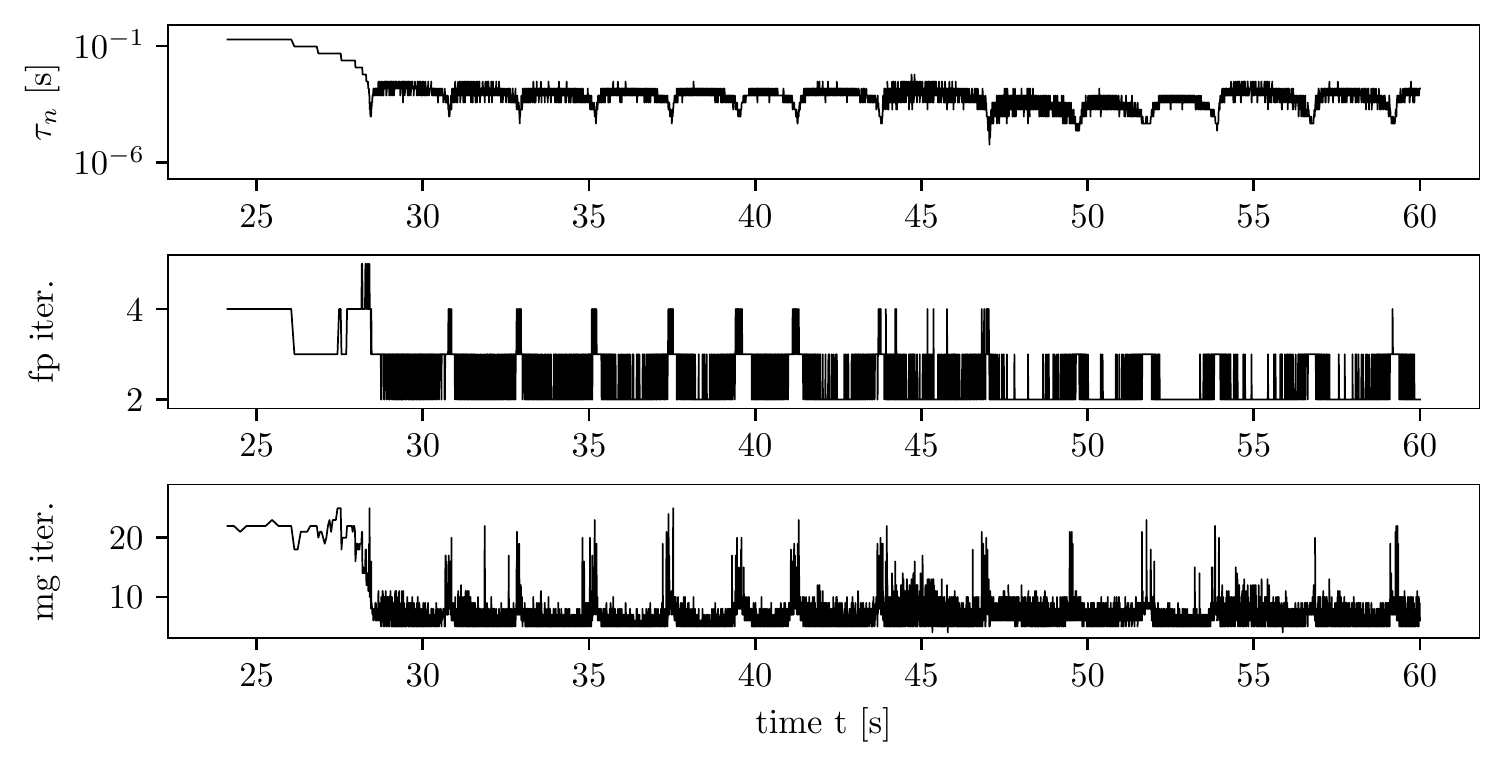}
    \vspace*{-25pt}
    \caption{Layered fault system:
    Adaptive time step selection and  performance of the algebraic solver.}
    \label{fig:MULLSOLVER}
\end{figure}

\subsubsection{Simulation results}
Figure~\ref{fig:MULLSLIPS} indicates quite interesting stress accumulation and release  along the different faults.
The four pictures show the mean value of the  relative velocity
on the faults $\frictionBoundary[4,5]$, $ \frictionBoundary[3,4]$, $\frictionBoundary[2,3]$,
and $\frictionBoundary[1,2]$ (top to bottom) over time.
In the first picture for the upper fault $\frictionBoundary[4,5]$,
we observe a sequence of almost periodic slip events with almost the same period and amplitude  as for the spring slider,
again after an initial loading phase of about  $\SI{26}{\second}$.
In both of the next two pictures, however, showing 
the average relative velocities
over the next two faults $ \frictionBoundary[3,4]$ and  $\frictionBoundary[2,3]$,
we see a highly oscillatory loading phase that seems to depict slip events several orders of magnitude smaller 
than on the top fault $\frictionBoundary[4,5]$ and might have saturated after a small jump at about $\SI{47}{\second}$ or lead to later slip events.
It is not clear at the moment whether the occurance and amplitude of these oscillations are physical or  due to numerical artifacts
which would motivate future numerical and experimental investigations.
As shown in the fourth picture, this jump of average relative velocity also occurs   
at the lowest fault  $\frictionBoundary[1,2]$, this time preceded by a rather stable loading phase.

Figure~\ref{fig:MULLLEVEL} shows the level lines of  approximate relative velocities 
along the upper fault  $\frictionBoundary[4,5]$ (horizontal axis)
evolving over six  time intervals (vertical axis) 
associated with the first six slip events (top left to bottom right).
All events exibit bilateral characteristics, i.e. ruptures nucleating towards the center of the fault and spreading towards both edges,
preceded by small foreshocks occurring in the middle of  $\frictionBoundary[4,5]$.
In most (but not all) cases the slip events are followed by small aftershocks.
These observations are in strong analogy with the results of the spring slider experiment as depicted in Figure~\ref{fig:SPRISLEVEL}.
\subsubsection{Adaptive time step selection and  performance of the algebraic solver}

The performance of adaptive time step selection and  of the algebraic solver 
as illustrated in Figure~\ref{fig:MULLSOLVER} hardly differs from the spring slider experiment.
Again, the slip events on $\frictionBoundary[4,5]$ are well-captured by adaptive time stepping 
that is reducing the time step by about 2 orders of magnitude.
The number of outer fixed point iterations  still ranges from 2 to 4 
and the sum of all inner multigrid iterations in each of these steps is bounded by 25, 
apart from slightly larger values at the end of the loading phase.
This strongly confirms the efficiency and robustness of our solution approach.

\section*{Statements and Declarations}

\subsection*{Funding}
This research has been funded by Deutsche Forschungsgemeinschaft (DFG)
through grant CRC 1114 ''Scaling Cascades in Complex Systems'', Project Number 235221301,
Project B01  ''Fault networks and scaling properties of deformation accumulation''.

\subsection*{Competing interests}
The authors have no competing interests to declare that are relevant to the
content of this article.

\subsection*{Data availability}
The datasets generated during and/or analysed during the current study
are available from the corresponding author on reasonable request.

\bibliographystyle{abbrv}
\bibliography{paper}

\bibliographystyle{abbrv}
\bibliography{paper}

\end{document}